\documentclass[11pt]{article}

\usepackage{amstext,amssymb,amsmath,amsbsy}
\usepackage{hyperref}
\usepackage{amscd}
\usepackage{amsfonts}
\usepackage{indentfirst}
\usepackage{verbatim}
\usepackage{amsmath}
\usepackage{amsthm}
\usepackage{enumerate}
\usepackage{graphicx}
\usepackage{color, soul}
\usepackage[OT1]{fontenc}
\usepackage[latin1]{inputenc}
\usepackage[english]{babel}
\usepackage{amssymb}
\usepackage{subfig}
\usepackage{algorithm}
\usepackage{algorithmic}

\usepackage{mathtools}
\DeclarePairedDelimiter\ceil{\lceil}{\rceil}

\newcommand{\R}{\mathbb{R}}

\newcommand{\ds}{\displaystyle}

\newcommand{\x}{{\bf x}}
\newcommand{\p}{{\bf p}}
\newcommand{\q}{{\bf q}}

\newcommand{\Div}{{\rm div}}

\newtheorem{Theorem}{Theorem}[section]
\newtheorem{Lemma}{Lemma}[section]

\newtheorem{Corollary}{Corollary}[section]
\newtheorem{remark}{Remark}[section]

\newtheorem*{Assumption*}{Assumption}
\newtheorem{Definition}{Definition}[section]
\newtheorem{problem}{Problem}[section]
\newtheorem*{problem*}{Problem}
\numberwithin{equation}{section}

\begin{document}

\title{Carleman estimates and the contraction principle for an inverse source problem for nonlinear hyperbolic equations}

\author{ Loc H. Nguyen\thanks{Department of Mathematics and Statistics, University of North Carolina at
Charlotte, Charlotte, NC, 28223, USA,  \texttt{loc.nguyen@uncc.edu} (corresponding
author).} \and Michael V. Klibanov\thanks{Department of Mathematics and Statistics, University of North Carolina at
Charlotte, Charlotte, NC, 28223, USA,  \texttt{mklibanv@uncc.edu.}}} 


\date{}
\maketitle
\begin{abstract}
	The main aim of this paper is to solve an inverse source problem for a general nonlinear hyperbolic equation. 
	Combining the  quasi-reversibility method and a suitable Carleman weight function, we define a map of which fixed point is the solution to the inverse problem.
	To find this fixed point, we define a recursive sequence with an arbitrary initial term by the same manner as in the classical proof of the contraction principle.
	Applying a Carleman estimate, we show that the sequence above converges to the desired solution with the exponential rate.
	Therefore, our new method  can be considered as an analog of the contraction principle.
	We rigorously study  the stability  of our method with respect to noise. 
	Numerical examples are presented.
\end{abstract}


\noindent{\it AMS subject classification: 	35R30;
65M32.

}

\section{Introduction}

Let $d \geq 2$ be the spatial dimension.
Let $F: \R^d \times [0, \infty) \times \R \times \R \times \R^d  \to \R$ be a function in the class $C^1$. 
For $u \in C^2({\R^d \times (0, \infty)})$, define the operator $\mathcal F$ as
\[
	\mathcal F u(\x, t) = F(\x, t, u, u_t, \nabla u) \quad \mbox{for all} \, (\x, t) \in \R^d \times (0, \infty).
\]
Let $c:\R^d \to [1, \overline c]$, for some given finite number $\overline c$, be a function in the class $C^2(\R^d)$.
Let $T > 0$ be a positive number representing the final time.
Consider the following initial value problem 
\begin{equation}
\left\{
	\begin{array}{ll}
	c(\x)u_{tt}(\x, t) = \Delta u(\x, t) + \mathcal F u(\x, t) 
	&(\x, t) \in \R^d \times (0, T),\\
	u_t(\x, 0) = 0 &\x \in \R^d,\\
	u(\x, 0) = \p(\x) &\x \in \R^d
	\end{array}
\right.
\label{main eqn}
\end{equation}
where $\p$ is the initial source function.
Let $\Omega$ be an open and bounded domain in $\R^d$ with smooth boundary $\partial \Omega$.
Define $\Omega_T = \Omega \times [0, T]$ and $\Gamma_T = \partial \Omega \times [0, T].$ 
Finding optimal conditions that guarantee the  existence and uniqueness results for the general problem like \eqref{main eqn} is challenging.
 Studying these properties for \eqref{main eqn} is out of the scope of this paper. We consider them as  assumptions.

	For the completeness, we provide here a set of conditions for $\p,$ $c$ and $F$ that guarantee the well-posedness of \eqref{main eqn}.
	Assume that $\p$ is smooth and has compact support and $c$ be a smooth and bounded function with $c \geq 1$. Assume that $F$ does not depend on the first derivatives of $u$
	 and that $|\mathcal F u| \leq C_1|u| + C_2$ for all functions $u \in H^2(\Omega \times (0, T))$ for some positive constants $C_1$ and $C_2$.
	Then, in this case, the  unique solvability for \eqref{main eqn} can be obtained by applying  \cite[Theorem 10.14]{Brezis:Springer2011}. This theorem is originally proved by Lions in \cite{Lions:sp1972}.
	We also refer the reader to the famous books \cite{Evans:PDEs2010, Ladyzhenskaya:sv1985} for the well-posedness of \eqref{main eqn} in the linear cases.

The main aim of this paper is to solve the following inverse problem:
\begin{problem}[An inverse source problem for nonlinear hyperbolic equations]
	 Determine the source function $\p(\x),$ $\x \in \Omega$, from the lateral Cauchy data
	\begin{equation}
		f(\x, t) = u(\x, t) \quad
		{\rm and}
		\quad
		g(\x, t) = \partial_\nu u(\x, t)
	\end{equation}  
	for all $(\x, t) \in \Gamma_T.$
	\label{ISP}
\end{problem}

Problem \ref{ISP} can be considered as the nonlinear version of the important problem arising from bio-medical imaging, called thermo/photo-acoustics tomography. 
The experiment leading to this problem is as follows, see \cite{Krugeretal:mp1995, Krugerelal:mp1999, Oraevskyelal:ps1994}.
One sends non-ionizing laser pulses or microwave to a biological tissue under inspection (for instance, woman's breast in mamography). A part of the energy will be absorbed and converted into heat, causing a thermal expansion and a subsequence ultrasonic wave propagating in space. The ultrasonic pressures on a surface around the tissue are measured.  Finding some initial information of the pressures from these measurements yields the structure inside this tissue.
The current techniques to solve the problem of thermo/photo-acoustics tomography are only for linear  hyperbolic equation.
We list here some widely used methods.
In the case when the waves propagate in the free space, one can find explicit reconstruction formulas in \cite{DoKunyansky:ip2018,  Haltmeier:cma2013, Natterer:ipi2012, Linh:ipi2009}, the time reversal method \cite{ KatsnelsonNguyen:aml2018, Hristova:ip2009, HristovaKuchmentLinh:ip2006, Stefanov:ip2009, Stefanov:ip2011}, the quasi-reversibility method \cite{ClasonKlibanov:sjsc2007, LeNguyenNguyenPowell:JOSC2021} and the iterative methods \cite{Huangetal:IEEE2013, Paltaufetal:ip2007, Paltaufetal:osa2002}.
The publications above study thermo/photo-acoustics tomography  for simple models for non-damping and isotropic media. 
The reader can find publications about thermo/photo-acoustics tomography for more complicated model involving a damping term or attenuation term \cite{Ammarielal:sp2012, Ammarietal:cm2011, Haltmeier:jmiv2019, Acosta:jde2018, Burgholzer:pspie2007, Homan:ipi2013, Kowar:SISI2014, Kowar:sp2012, Nachman1990}.
Up to the author's knowledge, this inverse source problem when the governing equation is general and nonlinear like \eqref{main eqn} is not solved yet.

Natural approaches to solve nonlinear inverse problem are based on optimization. 
These approaches are local in the sense that they only deliver solutions if good initial guesses of the solutions are given.
The local convergence does not always hold unless some additional conditions are imposed.
We refer the reader to \cite{Hankeetal:nm1995} for a condition that guarantees the local convergence of the optimization method using Landweber iteration.
The question how to compute the solutions to nonlinear problems without requesting a good initial guess is very interesting, challenging and significant in the scientific community.
There is a general framework to globally solve nonlinear inverse problem, called convexification. 
The main idea of the convexification method is to employ some suitable Carleman weight functions to convexify  the mismatch functional.
The convexified phenomenon is rigorously proved by employing the well-known Carleman estimates.
Several versions of the convexification method  \cite{KlibanovNik:ra2017, KhoaKlibanovLoc:SIAMImaging2020, Klibanov:sjma1997, Klibanov:nw1997, Klibanov:ip2015,  KlibanovKolesov:cma2019, KlibanovLiZhang:ip2019, KlibanovLeNguyenIPI2021,  KlibanovLiZhang:SIAM2019, LeNguyen:preprint2021} have
been developed since it was first introduced in \cite{KlibanovIoussoupova:SMA1995}. 
Especially, the convexification was successfully tested with experimental data in \cite{VoKlibanovNguyen:IP2020, Khoaelal:IPSE2021, KlibanovLeNguyenIPI2021} for the inverse scattering problem in the frequency domain given only back scattering data.
We consider the convexification method as the first generation of numerical method based on Carleman estimates to solve nonlinear inverse problems.
Although effective, the convexification method has a drawback. 
It is time consuming. 
We therefore propose a new method based on the fixed-point iteration, the contraction principle and a suitable Carleman estimate to globally solve nonlinear inverse problems.
This new method is considered as the second generation of globally convergent numerical methods to solve nonlinear inverse problems while the convexification method is the first one.

Our approach to solve Problem \ref{ISP} is to reformulate it to be of the form $\Phi(x) = x$ for some operator $\Phi$.
The general framework to find the fixed point of $\Phi$ in this paper is to construct a sequence by the same the sequence defined in the classical proof of the contraction principle.
Take any initial solution to the nonlinear inverse problem, denoted by $x_0$ and then recursively define $x_{n} = \Phi(x_{n-1})$ for $n > 0$. 
In general, the convergence of $\{x_n\}_{n \geq 0}$ to the fixed point of $\Phi$ is not guaranteed.
However, we figured out recently in \cite{LeNguyen:2020}, see also the strong publications \cite{BAUDOUIN:SIAMNumAna:2017, Baudouin:SIAM2021}, that if we include some suitable Carleman weight function in the definition of $\Phi$ then the  convergence holds true.
To the best knowledge to the authors, the idea of using Carleman estimate to construct a contraction map was first introduced in \cite{Baudouin2013}.
In the current paper, we define $\Phi$ by including the Carleman weight function in \cite[Lemma 6.1]{Klibanov:anm2015}, \cite[Theorem 1.10.2]{BeilinaKlibanovBook} and \cite[Section 2.5]{KlibanovLiBook}  into the classical quasi-reversibility method,
which was first introduced in \cite{LattesLions:e1969}.
The proof of the contractional behavior of $\Phi$ relies on the Carleman estimate in \cite[Lemma 6.1]{Klibanov:anm2015} and \cite[Theorem 1.10.2]{BeilinaKlibanovBook}.
The main theorem in this paper confirms that the recursive sequence above converges to the true solution. The stability  with respect to the noise contained in the given data is of the H\"older type.
The strengths of our new approach include the facts that
\begin{enumerate}
	\item it does not require a good initial guess;
	\item it is quite general in the sense that no special structure is imposed the nonlinearity $F$; 
	\item the rate of convergence is fast.
\end{enumerate}

The paper is organized as follows.
In Section \ref{sec 2}, we introduce our iteration.
In Section \ref{sec 3}, we prove the convergence of the sequence generated in Section \ref{sec 2}.
In Section \ref{num}, we present some numerical tests.
Section \ref{sec 5} is for concluding remarks.

\section{An iterative method to solve Problem \ref{ISP}} \label{sec 2}

Problem \ref{ISP} can be reduced to the problem of computing a function $u = u(\x, t)$ defined on $\overline {\Omega \times [0, T]}$ satisfying the following problem involving a nonlinear hyperbolic equation, the lateral Cauchy data and an initial condition
\begin{equation}
	\left\{
		\begin{array}{ll}
			c(\x) u_{tt}(\x, t) = \Delta u(\x, t) + \mathcal F u(\x, t) &(\x, t) \in \Omega_T,\\
			u(\x, t) = f(\x, t) & (\x, t) \in \Gamma_T,\\
			\partial_{\nu}u(\x, t) = g(\x, t) & (\x, t) \in \Gamma_T,\\
			u_t(\x, 0) = 0 & \x \in \Omega.
		\end{array}
	\right.
	\label{main eqn inverse}
\end{equation}
In the case of success, one can set $\p(\x) = u(\x, 0)$ for all $\x \in \Omega$ as the computed solution to Problem \ref{ISP}.
\begin{remark}
	A strength of the method in this paper is that we do not require any special condition on the nonlinearity $F$. In contrast, we assume that \eqref{main eqn inverse} has a unique and smooth solution.
	Since this paper focuses on computational method to solve inverse problem, studying the existence, uniqueness and regularity for \eqref{main eqn inverse} is out of the scope of this paper.
	We refer the reader to \cite{Brezis:Springer2011, Evans:PDEs2010} for these properties of solutions to hyperbolic equations.
\end{remark}
Due to the presence of the nonlinearity $\mathcal F$, numerically solving \eqref{main eqn inverse} is extremely challenging.
Conventional methods to do so are based on optimization.
That means, one minimizes some mismatch functionals, for e.g., 
\[
	u \mapsto \int_{\Omega_T}\Big|c(\x) u_{tt}(\x, t) - \Delta u(\x, t) - \mathcal F u(\x, t)\Big|^2 + \mbox{a regularization term}
\]
subject to the boundary and initial constraints in the last three lines of \eqref{main eqn inverse}.
The computed solution to \eqref{main eqn inverse} is set as the obtained minimizer.
However, due to the nonlinearity $\mathcal F$, such a functional is not convex.
It might have multiple minima and ravines.
Therefore, the success of the methods based on optimization depends on the good initial guess. 
In reality, good initial guesses are not always available.
Motivated by this challenge, we propose to construct a fixed-point like iterative sequence based on the quasi-reversibility method and a Carleman estimate in \cite[Lemma 6.1]{Klibanov:anm2015} and \cite[Theorem 1.10.2]{BeilinaKlibanovBook}. 
The convergence of this sequence to compute solution to \eqref{main eqn inverse} will be proved by the contraction principle and a Carleman estimate.
Our Carleman-based iterative algorithm to numerically solve \eqref{main eqn inverse} is as follows.

Let $\x_0$ be a point in $\R^d\setminus \overline \Omega.$ 
For $\lambda > 0$ and $\eta > 0$, define the Carleman weight function
\begin{equation}
	W_{\lambda, \eta}(\x, t) = e^{\lambda(|\x - \x_0|^2 - \eta t^2)} \quad
	 \mbox{for all } (\x, t) \in \Omega_T.
	 \label{Car w}
\end{equation}
Let $p > \ceil{(d + 1)/2} + 2$ where $\ceil{(d + 1)/2}$ is the smallest integer that is greater than or equal to $(d + 1)/2$.
Then, $H^p(\Omega_T)$ is continuously embedded into $C^2(\Omega_T).$
The convergence of our method to solve \eqref{main eqn inverse}  is based on the following Carleman weighted norms.
\begin{Definition}[Carleman weighted Sobolev norms]
	For all $\lambda \in (1, \infty), \eta \in [0, 1)$ and $s \in \mathbb N,$ we define the norm
	\begin{equation}
	\|v\|_{H_{\lambda, \eta}^s(\Omega_T)} = \Big[
		\sum_{|\alpha| = 0}^s\int_{\Omega_T} W^2_{\lambda, \eta}(\x, t)|D_\alpha v|^2dt d\x
	\Big]^{1/2} 
	\quad 
	\mbox{for all } v \in H^s(\Omega_T).
	\label{Hsnorm}
\end{equation}
The notation $D_\alpha$ in \eqref{Hsnorm} is understood in the usual sense. 
That means, $\alpha = (\alpha_1, \dots, \alpha_d, \alpha_{d + 1})$ being a $d + 1$ dimensional vector of nonnegative integers, $|\alpha| = \ds\sum_{i = 1}^{d + 1}\alpha_i$ and $D_\alpha = \ds\frac{\partial^{|\alpha|}}{\partial_{x_1}^{\alpha_1} \dots \partial_{x_d}^{\alpha_d} \partial_t^{\alpha_{d + 1}}}.$ 
\label{def1}
\end{Definition}

Assume that the set of admissible solutions
\[
	H = \{u \in H^p(\Omega_T): u_t|_{\Omega \times \{0\}} = 0, u|_{\Gamma_T} = f \,{\rm and} \, \partial_{\nu} u|_{\Gamma_T} = g\}
\]
is nonempty. 
We construct a sequence $\{u_n\}_{n \geq 0} \subset H$ that converges to the solution to \eqref{main eqn inverse}.
Take an arbitrary function $u_0 \in H$. 
Assume by induction that $u_n$ is known.
We set $u_{n+1} = \Phi(u_n)$
where $\Phi(u_n)$ is the unique minimizer of $J_n: H \to \R$ defined as
\begin{equation}
	J_n(v) = \int_0^T\int_{\Omega} W_{\lambda, \eta}^2(\x, t) |c(\x)v_{tt} - \Delta v - \mathcal F(u_n)|^2d\x dt
	+ \epsilon \|v\|_{H^p_{\lambda, 0}(\Omega_T)}^2 
	\quad 
	\mbox{for all } \, v \in H.
	\label{Jn}
\end{equation}
See \eqref{Hsnorm} for the norm of the regularity term $\epsilon \|v\|_{H^p_{\lambda, 0}(\Omega_T)}^2$.
Due to the compact embedding from $H^p(\Omega_T)$ to $H^2(\Omega_T),$ the functional $J_n$ is weakly lower semi-continuous. 
The presence of the regularization term guarantees that $J_n$ is coercive. 
Hence, $J_n$ has a minimizer $u_{n+1}$ on the close and convex set $H$ in $H^p(\Omega_T)$.
The uniqueness of $u_{n+1}$ is due to the strict convexity of $J_n$.
An alternative way to obtain the existence and uniqueness of the minimizer is to argue similarly to the proofs of Theorem 2.1 in \cite{KlibanovNik:ra2017} or Theorem 4.1 in \cite{KlibanovLeNguyenIPI2021}.
We will prove that the sequence $\{u_n\}_{n \geq 0}$ converges to the solution to \eqref{main eqn inverse}.
This suggests Algorithm \ref{alg} to solve Problem \ref{ISP}.
\begin{algorithm}[h!]
\caption{\label{alg} A procedure to compute the numerical solution to Problem \ref{ISP}}
	\begin{algorithmic}[1]
	\STATE \label{s1} 
	Choose a number $n_0$ that represents the maximal number of iteration. Choose an arbitrary initial solution $u_0 \in H.$
	\FOR{$n = 0$ to $n_0$} \label{n0}
		\STATE \label{step update} Minimize the function $J_n$, defined in \eqref{Jn}, on $H$. 
		\STATE Set $u_{n+1}$ as the obtained minimizer.
	\ENDFOR
	\STATE Set the computed solution to \eqref{main eqn inverse} as $u_{\rm comp} = u_{n_0}.$
	\STATE Set the computed solution to Problem \ref{ISP} as $\p_{\rm comp}(\x) = u_{\rm comp}(\x, 0)$ for all $\x \in \Omega.$ 
\end{algorithmic}
\end{algorithm}

\begin{remark}
 The presence of the Carleman weight functions $W_{\lambda, \eta}(\x, t)$ and that of the weight function $W_{\lambda, 0}(\x, t)$ in the regularization term  in \eqref{Jn} play the key role in proving the convergence of the sequence $\{u_n\}_{n \geq 0}$ to the true solution to $\eqref{main eqn}.$ See Theorem \ref{thm} for details. 
\end{remark}

\begin{remark}
	One can replace the homogeneous initial condition for $u_t(\x, 0)$ in \eqref{main eqn}  by the non-homogenous one $u_t(\x, 0) = \q(\x)$ for all $\x \in \R^d$ for some given function $\q$. 
	The analysis leading to the numerical method and the proof for the convergence of the method will be the same.
\end{remark}

\begin{remark}
	Besides Algorithm \ref{alg}, one can apply the convexification method, first introduced in \cite{KlibanovIoussoupova:SMA1995}, to globally solve \eqref{main eqn inverse} and Problem \ref{ISP}. 
	In fact, by  combining the arguments in \cite{KlibanovNik:ra2017, ClasonKlibanov:sjsc2007, KlibanovLeNguyenIPI2021}, one can prove that:
	\begin{enumerate}
	\item for any large ball $B(M) = \{u \in H^p(\Omega_T): \|u\|_{H^s(\Omega_T)} < M\}$ of $H^s(\Omega_T)$ where $M$ is an arbitrarily large number,   the functional
	\[
	u \mapsto \int_{\Omega_T} W^2_{\lambda, \eta}(\x, t)\Big|c(\x) u_{tt}(\x, t) - \Delta u(\x, t) - \mathcal F u(\x, t)\Big|^2 + \mbox{ a regularization term}
	\]
	is strictly convex in the large ball of $H^s(\Omega_T)$ for all $\lambda$ sufficiently large and $\eta$ sufficiently small;
	\item the unique minimizer of this functional in $B(M)$ is a good approximation of the desired solution.
	\end{enumerate}
	In theory, both Algorithm \ref{alg} and the convexification method deliver reliable solution to \eqref{main eqn inverse} without requesting initial guess.	 
	It is interesting to verify the efficiency of the convexification method while that of Algorithm \ref{alg} is confirmed in this paper. 
	Developing the convexification method for \eqref{main eqn inverse} and verify it numerically serve as one of our future research topics.
\end{remark}

\section{The convergence theorem} \label{sec 3}

In this section, we employ a Carleman estimate, see \cite[Lemma 6.1]{Klibanov:anm2015} and \cite[Theorem 1.10.2]{BeilinaKlibanovBook} to prove the convergence of the sequence $\{u_n\}_{n \geq 0}$ obtained by Algorithm \ref{alg} to the true solution to \eqref{main eqn inverse}.

\subsection{The statement of the convergence theorem}

Since $\Omega$ is a bounded domain, it is contained in a large disk $D = B(R) = \{\x \in \R^d: |\x| < R\}$ for some positive number $R$. 
Before stating the convergence theorems,
we  recall an important Carleman estimate.
For $\eta > 0$ and $\varepsilon > 0$, define
\begin{equation}
	D_{\eta, \varepsilon} = \{(\x, t) \in D \times [0, \infty): |\x - \x_0|^2 - \eta t^2 > \varepsilon\}.
\end{equation}
The following lemma plays a key role in our analysis.
\begin{Lemma}[Carleman estimate]
	Assume that the point $\x_0$ in the definition of the Carleman weight function in \eqref{Car w} is in the set $D \setminus \overline \Omega$
	and that the coefficient $c$ in \eqref{main eqn} satisfies the condtion
\begin{equation}
	\langle \x - \x_0, \nabla c(\x) \rangle \geq 0  \quad \mbox{for all } \x \in D. 
	\label{reg c}
\end{equation}
	Then, there exists a sufficiently small number $\eta_0 = \eta_0(\x_0, D, \|c\|_{C^1(\overline D)}) \in (0, 1)$ such that for any $\eta \in (0, \eta_0]$, one can choose a sufficiently large number $\lambda_0 = \lambda_0(D, \eta, c, \x_0) > 1$ and a number $C = C(D, \eta, c, \x_0)$ such that for all $z \in C^2(\overline{D_{\eta, \varepsilon}})$  and for all $\lambda \geq \lambda_0$, the following pointwise Carleman estimate holds true
	\begin{multline}
		W_{\lambda, \eta}^2(\x, t)|c(\x) z_{tt}(\x, t) - \Delta z(\x, t)|^2  
		\geq 
		C \lambda W_{\lambda, \eta}^2(\x, t) \Big(
			|\nabla z(\x, t)|^2 + |z_t(\x, t)|^2 + \lambda^2 |z(\x, t)|^2
		\Big)
		\\
		+ \Div(Z(\x, t)) + Y_t(\x, t) 
	\label{Car est}
	\end{multline}
	for $(\x, t) \in D_{\eta, \varepsilon}$. 
	The vector valued function $Z$ and the function $Y$ satisfy
	\begin{equation}
		|Z(\x, t)| \leq C \lambda^3 W_{\lambda, \eta}^2(\x, t)(|\nabla z(\x, t)|^2 + |z_t(\x, t)|^2 + |z(\x, t)|^2)
	\label{Z}
	\end{equation}
	and
	\begin{multline}
		|Y(\x, t)| \leq C \lambda^3 W_{\lambda, \eta}^2(\x, t)\Big[|t|(|\nabla z(\x, t)|^2 + |z_t(\x, t)|^2 + |z(\x, t)|^2) 
		\\
		+ \big(|\nabla z(\x, t)| + |z(\x, t)|\big) |z_t(\x, t)|\Big].
	\label{Y}
	\end{multline}
	In particular, if either $z(\x, 0) = 0$ or $z_t(\x, 0) = 0,$ then $Y(\x, 0) = 0$.
\label{lem Car}
\end{Lemma}
We refer the reader to \cite[Lemma 6.1]{Klibanov:anm2015}, \cite[Theorem 1.10.2]{BeilinaKlibanovBook} and \cite[Section 2.5]{KlibanovLiBook} for the proof of Lemma \ref{lem Car}. 
The proof of this Carleman weight function in the case $c \equiv 1$ is given in \cite{Lavrentiev:AMS1986}.
In those publications, the Carleman estimate \eqref{Car est} holds true for the function $z$ in $C^2(\overline {D_{\eta, \varepsilon}^{\pm}})$ where
\[
	D_{\eta, \varepsilon}^{\pm} = \{
		(\x, t) \in D \times [-T, T]: |\x - \x_0|^2 - \eta t^2 > \varepsilon
	\}.
\]
It is possible to bring in the exact formulas of two important functions $Y$ and $Z$ from the proof of Theorem 1.10.2 in \cite{BeilinaKlibanovBook}. However, we do not do so since the presence of these formulas has no contribution to the proof of our convergence theorem. Only their estimates in \eqref{Z} and \eqref{Y} are useful.
Another reason for us to not write those formulas is that they are long. 
For example, $Y$ is the sum of many complicated functions; for example, one of which is
\[
	4\lambda c(\x)^2 \eta t |v_t(\x, t)|^2 - 16 \lambda c(\x) \eta\Big(
		t|\x - \x_0| - c\eta^2 t^2 + tO(\lambda^{-1})
	\Big) v^2(\x, t)
\]
where $v(\x, t) = u(\x, t) W_{\lambda, \eta}^2(\x, t)$ for $(\x, t) \in Q_T$.
Writing some lengthy formulas having no contribution makes the presentation of the paper less concentrated. 
So, rather than writing these lengthy formulas, we focus on their important bounds in \eqref{Z} and \eqref{Y}.
The second author also mentioned that the  occurrence of long formula is typical in derivation of Carleman estimates, see  \cite[Sections 2.3 and Section 2.5]{KlibanovLiBook}.

Let $\eta_0$ be the number in Lemma \ref{lem Car}.
Fix $\varepsilon > 0$ and $\eta \in (0, \eta_0]$.
Since $\x_0 \not \in \Omega$, we can choose both of these two numbers sufficiently small;
say for e.g., 
$
	\varepsilon < \frac{P_1^2}{2}
$
and
$ \eta < \frac{P_1^2}{2T^2}
$
where $P_1 = \min_{\x \in \overline \Omega_0}\{|\x - \x_0|\}$, such that 
\begin{equation}	
	|\x - \x_0|^2 - \eta t^2 > \varepsilon \quad
	\mbox{for all } (\x, t) \in \Omega_T.
	\label{3.7777}
\end{equation}

\begin{remark}
With the choice of $\eta$ and $\varepsilon$ above, the Carleman estimate \eqref{Car est} holds true for all $(\x, t) \in \Omega_T.$
\label{rm 3.1}
\end{remark}


We now state the  main theorem of the paper.
Let $\delta > 0$ be a noise level. By noise level $\delta$, we assume that there is a function $\mathcal E$ in $H^p(\Omega_T)$ satisfying
\begin{equation}
	\|\mathcal E\|_{H^p(\Omega_T)} \leq \delta, 
	\quad
	f|_{\Gamma_T} = f^* +  \mathcal E|_{\Gamma_T},
	\quad
	g|_{\Gamma_T} = g^* +  \partial_{\nu} \mathcal E|_{\Gamma_T}
	\label{error function}
\end{equation}
where $f^*$ and $g^*$ are the noiseless version of the data $f$ and $g$ respectively.
The true solution to \eqref{main eqn inverse}, with $f$ and $g$ replaced by $f^*$ and $g^*$ respectively, is denoted by $u^*$. That means, $u^*|_{\Gamma_T} = f^*$ and $\partial_{\nu} u^*|_{\Gamma_T} = g^*.$

We first establish the convergence analysis when $\|F\|_{C^1} < \infty$. The case $\|F\|_{C^1} = \infty$ will be considered later by  truncating $F$ using a cut off function as in \eqref{cutoff}.

\begin{Theorem}[The convergence of $\{u_{n}\}_{n \geq 0}$ to $u^*$ as the noise tends to $0$]
Assume that $\|F\|_{C^1} $ is a finite number  and that the function $c$ satisfies \eqref{reg c}.
	Fix $\varepsilon \in (0, 1)$ and let $\eta \in (0, \eta_0]$ be such that \eqref{3.7777} holds true. 
	Fix $\gamma \in (0, 1).$
	Then, there exist $\lambda_0 > 1$ depending only on $D$, $\eta$, $c$, $\x_0$, $\varepsilon$ and $\epsilon$ such that for all $\lambda > \lambda_0$, 
we have
\begin{equation}
	\|u_{n + 1} - u^*\|^2_{H^1_{\lambda, \eta}(\Omega_T)}  \leq \theta^{n+1} \|u_{0} - u^*\|^2_{H^1_{\lambda, \eta}(\Omega_T)} + C(\delta^{2 - 2\gamma} + \epsilon \|u^*\|_{H^p_{\lambda, 0}(\Omega_T)}^2)
	\label{conv}
\end{equation}
provided that  $\delta < e^{-\lambda P_2^2/\gamma}$
	where $P_2 = \max_{\x \in \overline\Omega}{|\x - \x_0|}.$ 
	Here, $C$ is a positive number depending only on $\Omega,$ $D$, $c$, $\x_0$ and $T$ and $\theta = C/\lambda \in (0, 1)$. 
	\label{thm}
\end{Theorem}

\begin{Corollary}
	By the trace theory, it follows from \eqref{conv} that for all $n \geq 0,$
	\begin{equation}
	\|p_{n + 1} - p^*\|^2_{L^2_{\lambda}(\Omega)}  \leq C_1\theta^{n+1} \|u_{0} - u^*\|^2_{H^1_{\lambda, \eta}(\Omega_T)} + C_2(\delta^{2 - 2\gamma} + \epsilon \|u^*\|_{H^p}(\Omega_T)^2)
	\label{p}
\end{equation}
where $C_1$ and $C_2$ are two positive numbers depending only on  $\Omega,$ $D$, $c$, $\x_0$ and $T$.
In \eqref{p}, $p_{n + 1}(\x) = u_{n + 1}(\x, 0)$, $p^*$ is the true source function, and
\[
	\|p_{n + 1} - p^*\|_{L^2_{\lambda}(\Omega)}
	= \Big[
		\int_{\Omega} e^{2 \lambda |\x - \x_0|^2} |p_{n + 1} - p^*|^2 d\x
	\Big]^{\frac{1}{2}}.
\]
\end{Corollary}

\begin{remark}
	In the convergence theorem above, we impose 
the condition that the true solution of \eqref{main eqn} and \eqref{main eqn inverse} belongs to the class $H^p(\Omega_T)$, 
	which can be embedded into the class $C^2(\overline \Omega).$ Therefore, the data $f(\x, t)$ and $g(\x, t)$ must be the restrictions of a smooth function and its normal derivative respectively on $\partial \Omega \times [0, T]$.
	By standard results in regularity, this assumption holds true if the initial condition $\p,$ the nonlinearity $F$ and the coefficient $c$ are sufficiently smooth.
	These assumptions might not be practical.
	However, solving inverse problem is usually challenging.
	 As a rule, the smoothness requirements of coefficients and source terms in the governing PDEs are not of a primary concern 
	 see, e. g. Theorem 4.1 in \cite{Romanov:VNU1986}  as well as papers \cite{Novikov:jfa1992, NovikovImrs2005}.
	On the other hand, we also have to impose a strong condition about the smoothness of the noise, see \eqref{error function}. 
	Those smoothness assumptions and noise regularization can be relaxed in the numerical study.
	The numerical results due to our new method are out of expectation even when the function $F$ is not smooth and the initial condition $\p$ is not continuous. 
	Moreover, in our numerical study, the noise is the function taking uniformly distributed random numbers, see \eqref{noisenoise}. 	
	That means the procedure in Algorithm \ref{alg} is stronger than what we can prove. 
	\label{rem3.3333}
\end{remark}

	In practice, the condition that $\|F\|_{C^1(\overline\Omega \times \R \times \R^d)} < \infty$ might not hold true.
	In this case, we need to know in advance a number $M$ such that the true solution $u^*$ to \eqref{main eqn inverse} belongs to the ball $B_M = \{v \in H^p(\Omega_T): \|v\|_{H^p(\Omega_T)} < M\}$. 
		This requirement does not weaken our result since $M$ can be arbitrarily large.
	In order to compute $u^*$, we proceed as follows.	
	Due to the Sobolev embedding theorem, $\|u^*\|_{C^1(\overline \Omega_T)} < C_3 M$ for some constant $C_3$.
	Let $\chi: \overline \Omega \times [0, \infty) \times \R \times \R \times \R^d \to \R$ be a function in the class $C^1$ and satisfy
	\begin{equation}
		\chi(\x, t, s_1, s_2, \p) = 
		\left\{
			\begin{array}{ll}
				1 & (s_1^2 + s^2 + |\p|^2)^{1/2}  \leq C_3 M,\\
				\in (0, 1) & C_3 M < (s_1^2 + s^2 + |\p|^2)^{1/2} < 2C_3 M,\\
				0 & (s_1^2 + s^2 + |\p|^2)^{1/2}  \geq 2C_3 M.
			\end{array}
		\right.	
		\label{cutoff}
	\end{equation}
	Then, it is obvious that $u^*$ solves
		\begin{equation}
	\left\{
		\begin{array}{ll}
			c(\x) u_{tt}(\x, t) = \Delta u(\x, t) +\overline {\mathcal F} u(\x, t) &(\x, t) \in \Omega_T,\\
			u(\x, t) = f(\x, t) & (\x, t) \in \Gamma_T,\\
			\partial_{\nu}u(\x, t) = g(\x, t) & (\x, t) \in \Gamma_T,\\
			u_t(\x, 0) = 0 & \x \in \Omega
		\end{array}
	\right.
	\label{3.9999}
\end{equation}
where
\[
	\overline {\mathcal F} u = \chi(\x, t, u, u_t, \nabla u) F(\x, t, u, u_t, \nabla u).
\]
Since the function $ \chi F$ is smooth and has compact support, it has a bounded $C^1$ norm.
We can compute $u^*$ by solving \eqref{3.9999} using Algorithm \ref{alg} and Theorem \ref{thm} for $\overline {\mathcal F}$.


We present the proof of  Theorem \ref{thm} in the next subsection.
The proof is motivated by \cite{BAUDOUIN:SIAMNumAna:2017, Baudouin:SIAM2021, LeNguyen:2020}.

\subsection{The proof of Theorem \ref{thm}}

Throughout the proof, $C$ is a generic constant depending only on $\Omega,$ $D$, $c$, $\x_0$, $T$, and $\|F\|_{C^1}$. 
In the proof, we use the Carleman weighted Sobolev norms $\|\cdot\|_{H^s_{\lambda, \eta}(\Omega_T)}$ defined in Definition \ref{def1} with some different values of $s$, $\lambda$ and $\eta$.
We split the proof into several steps.

{\bf Step 1.}
Define 
\[
	H_0 = \{u \in H^p(\Omega_T): u_t|_{\Omega \times \{0\}} = 0, u|_{\Gamma_T} = 0, \mbox{and}\, \partial_\nu u|_{\Gamma_T} = 0\}.
\]
Fix $n > 0$.
Since $u_{n + 1}$ is the minimizer of $J_n$, defined in \eqref{Jn},  on $H$, by the variational principle, we have for all $h \in H_0$
\begin{multline}
\int_{\Omega_T}	W_{\lambda, \eta}^2(\x, t) \big(c(\x)(u_{n+1})_{tt} - \Delta u_{n + 1} - \mathcal F(u_n)\big)\big(c(\x)h_{tt} - \Delta h \big)d\x dt
	\\
	+ \epsilon \langle u_{n+1}, h\rangle_{H^p_{\lambda, 0}(\Omega_T)} = 0.
	\label{3.6}
\end{multline}
On the other hand, since $u^*$ is the true solution to \eqref{main eqn inverse}, for all $h \in H_0$, we have
\begin{multline}	
\int_{\Omega_T}	W_{\lambda, \eta}^2(\x, t) \big(c(\x)u^*_{tt} - \Delta u^* - \mathcal F(u^*)\big)\big(c(\x)h_{tt} - \Delta h \big)d\x dt
	+ \epsilon \langle u^*, h\rangle_{H^p_{\lambda, 0}(\Omega_T)}
	\\
 = \epsilon \langle u^*, h\rangle_{H^p_{\lambda, 0}(\Omega_T)}.
 \label{3.7}
\end{multline}
Combining \eqref{3.6} and \eqref{3.7}, for all $h \in H_0$, we have
\begin{multline}
 \int_{\Omega_T}	W_{\lambda, \eta}^2(\x, t) \big(c(\x)(u_{n+1} - u^*)_{tt} - \Delta (u_{n+1} - u^*) - (\mathcal F(u_n) - \mathcal F(u^*))\big)\big(c(\x)h_{tt} - \Delta h \big)d\x dt
	\\
	+ \epsilon \langle u_{n+1} - u^*, h\rangle_{H^p_{\lambda, 0}(\Omega_T)}
	=  -\epsilon \langle  u^*, h\rangle_{H^p_{\lambda, 0}(\Omega_T)}.
	\label{3.8}
\end{multline}

Due to \eqref{error function}, the function 
\begin{equation}
	z = u_{n+1} - u^* - \mathcal{E} 
	\quad \mbox{or} 
	\quad z + \mathcal{E} = u_{n + 1} - u^*
	\label{test fn h}
\end{equation} is in $H_0$.
Using this test function $h = z$ for the identity \eqref{3.8}, we have
\begin{multline}
 \int_{\Omega_T}	W_{\lambda, \eta}^2(\x, t) \big(c(\x)(z + \mathcal E)_{tt} - \Delta (z + \mathcal E) - (\mathcal F(u_n) - \mathcal F(u^*))\big)
	\\
	\big(c(\x)z_{tt} - \Delta z \big)d\x dt
	+ \epsilon \langle z + \mathcal E, z\rangle_{H^p_{\lambda, 0}(\Omega_T)}
	=  -\epsilon \langle  u^*, z\rangle_{H^p_{\lambda, 0}(\Omega_T)}.
	\label{3.16}
\end{multline}
Thus, by \eqref{3.16}, we have
\begin{multline}
 \int_{\Omega_T}	W_{\lambda, \eta}^2(\x, t) \big|c(\x)z_{tt} - \Delta z  \big|^2 d\x dt
	+ \epsilon \|z\|^2_{H^p_{\lambda, 0}(\Omega_T)}
	\\
	 = 
 \int_{\Omega_T}	W_{\lambda, \eta}^2(\x, t) \big(\mathcal F(u_n) - \mathcal F(u^*)\big)\big(c(\x)z_{tt} - \Delta z \big) d\x dt
	\\
	-  \int_{\Omega_T}	W_{\lambda, \eta}^2(\x, t) \big(c(\x)\mathcal E_{tt} - \Delta \mathcal E \big)\big(c(\x)z_{tt} - \Delta z \big)
	 - \epsilon \langle u^*, z\rangle_{H^p_{\lambda, 0}(\Omega_T)}.
	\label{3.9}
\end{multline}
Using the inequality $|ab| \leq 4a^2 + b^2/16$, we deduce from \eqref{3.9} that
\begin{multline}
	 \int_{\Omega_T}	W_{\lambda, \eta}^2(\x, t) \big|c(\x)z_{tt} - \Delta z  \big|^2 d\x dt
	+ \epsilon \|z\|^2_{H^p_{\lambda, 0}(\Omega_T)}
	 \leq 
	C\Big[  \int_{\Omega_T}	W_{\lambda, \eta}^2(\x, t) \big|\mathcal F(u_n) - \mathcal F(u^*)\big|^2 d\x dt
	\\
	+ \int_{\Omega_T}	W_{\lambda, \eta}^2(\x, t) \big|c(\x)\mathcal E_{tt} - \Delta \mathcal E \big|^2d\x dt
	 +\epsilon\|u^*\|^2_{H^p_{\lambda, 0}(\Omega_T)}\Big].
	\label{3.10}
\end{multline}

Since $\|F\|_{C^1(\overline\Omega \times \R \times \R^d)} < \infty$, we can find a constant $C$ such that
\[
	|\mathcal F(u_n) - \mathcal F(u^*)|^2 \leq C [|u_n - u^*|^2 + |(u_n - u^*)_t|^2 + |\nabla (u_n - u^*)|^2] 
	\quad \mbox{in } \Omega_T.
\]
Using this and \eqref{3.10}, we have
\begin{multline}
	 \int_{\Omega_T}	W_{\lambda, \eta}^2(\x, t) \big|c(\x)z_{tt} - \Delta z  \big|^2 d\x dt
	+ \epsilon \|z\|^2_{H^p_{\lambda, 0}(\Omega_T)}
	\\
	 \leq 
	C\Big[  \int_{\Omega_T}	W_{\lambda, \eta}^2(\x, t) \big( |u_{n} - u^*|^2  + |(u_n - u^*)_t|^2 + |\nabla (u_n - u^*)|^2\big)  d\x dt
	\\
	+  \int_{\Omega_T}	W_{\lambda, \eta}^2(\x, t) \big|c(\x)\mathcal E_{tt} - \Delta \mathcal E \big|^2d\x dt
	 +\epsilon\|u^*\|^2_{H^p_{\lambda, 0}(\Omega_T)}\Big].
	\label{3.11}
\end{multline}

{\bf Step 2.} 
Due to Remark \ref{rm 3.1}, the Carleman estimate \eqref{Car est} holds true for the function $z$ in the whole set $\Omega_T$. 
Integrating  \eqref{Car est} on $\Omega_T$, we have
\begin{multline}
		\int_{\Omega_T}W_{\lambda, \eta}^2(\x, t)|c(\x) z_{tt}(\x, t) - \Delta z(\x, t)|^2  dt d\x
		\\
		\geq 
		C \lambda \int_{\Omega_T} W_{\lambda, \eta}^2(\x, t) \Big(
			|\nabla z(\x, t)|^2 + |z_t(\x, t)|^2 + \lambda^2 |z(\x, t)|^2
		\Big) dtd\x
		\\
		+ \int_{\Omega_T}\Div(Z(\x, t)) dtd\x 
		+ \int_{\Omega_T}Y_t(\x, t)  dtd\x.
		\label{3.15}
\end{multline}
Since $z \in H_0,$ $Z(\x, t) = 0$ for all $\x \in \partial \Omega$.
 By the divergence theorem, we have 
\begin{equation}
	\int_{\Omega_T}\Div(Z(\x, t)) dtd\x = 0.
	\label{3.17}
\end{equation}
Using the Fundamental theorem of Calculus, we have
\begin{equation}
	\int_{\Omega_T}Y_t(\x, t)  dtd\x 
	= \int_{\Omega} \int_0^T Y_t(\x, t)  dtd\x  
	= \int_{\Omega} [Y(\x, T) - Y(\x, 0)]d\x. 
	\label{3.18}
\end{equation}
Since $z_t(\x, 0) = 0$ for all $\x \in \Omega$, we have $Y(\x, 0) = 0$. 
Using \eqref{Y}, we deduce from \eqref{3.18} that
\begin{equation}
	\Big|
		\int_{\Omega_T}Y_t(\x, t)  dtd\x
	\Big|
	\leq   C \lambda^3
	\int_{\Omega}  W_{\lambda, \eta}^2(\x, T)\big(|\nabla z(\x, T)|^2 + |z_t(\x, T)|^2 + |z(\x, T)|^2 
		\Big)d\x.
	\label{3.1818}		
\end{equation}
It follows from \eqref{3.15} and \eqref{3.1818} that
\begin{multline}
	\int_{\Omega_T}W_{\lambda, \eta}^2(\x, t)|c(\x) z_{tt}(\x, t) - \Delta z(\x, t)|^2  dt d\x
		\\
		\geq 
		C \lambda \int_{\Omega_T} W_{\lambda, \eta}^2(\x, t) \Big(
			|\nabla z(\x, t)|^2 + |z_t(\x, t)|^2 + \lambda^2 |z(\x, t)|^2
		\Big) dtd\x
		\\
		-  C \lambda^3
	\int_{\Omega}  W_{\lambda, \eta}^2(\x, T)\big(|\nabla z(\x, T)|^2 + |z_t(\x, T)|^2 + |z(\x, T)|^2 
		\Big)d\x.
		\label{3.19}
\end{multline}
Using the trace theory on $H^p(\Omega_T) \hookrightarrow H^1(\Omega \times \{t = T\}),$ we can find a constant $C_1$ such that
\begin{align}
	\lambda^3 e^{-\lambda \eta T^2} \|z\|^2_{H^p_{\lambda, 0}(\Omega_T)} &=
	 	\lambda^3 e^{-\lambda \eta T^2}\int_{\Omega_T}\sum_{|\alpha| = 0}^p e^{2\lambda|\x - \x_0|^2} |\partial_{\alpha} z|^2 dtd\x \notag
	 \\
	&\geq C_1 \lambda^3  \int_{\Omega} W_{\lambda, \eta}^2(\x, T) \Big(
		|\nabla z(\x, T)|^2 + |z_t(\x, T)|^2 + |z(\x, T)|^2
	\Big)dt d\x
	\label{3.20}
\end{align}
Combining \eqref{3.19} and \eqref{3.20} gives 
\begin{multline}
	\int_{\Omega_T}W_{\lambda, \eta}^2(\x, t)|c(\x) z_{tt}(\x, t) - \Delta z(\x, t)|^2  dt d\x 
	+ \varepsilon \|z\|^2_{H^p_{\lambda, 0}(\Omega_T)}
		\\
		\geq 
		C \lambda \int_{\Omega_T} W_{\lambda, \eta}^2(\x, t) \Big(
			|\nabla z(\x, t)|^2 + |z_t(\x, t)|^2 + \lambda^2 |z(\x, t)|^2
		\Big) dtd\x
		\\
		+ \varepsilon \|z\|^2_{H^p_{\lambda, 0}(\Omega_T)} -  C  \lambda^3 e^{-\lambda \eta T^2} \|z\|^2_{H^p_{\lambda, 0}(\Omega_T)}.
		\label{3.21}
\end{multline}

{\bf Step 3.}
It follows from \eqref{3.11} and \eqref{3.21} that
\begin{multline}
	C \lambda \int_{\Omega_T} W_{\lambda, \eta}^2(\x, t) \Big(
			|\nabla z(\x, t)|^2 + |z_t(\x, t)|^2 + \lambda^2 |z(\x, t)|^2
		\Big) dtd\x
		+ \varepsilon \|z\|^2_{H^p_{\lambda, 0}(\Omega_T)} 
		\\
		-  C  \lambda^3 e^{-\lambda \eta T^2} \|z\|^2_{H^p_{\lambda, 0}(\Omega_T)}
	\leq 
	 \int_{\Omega_T}	W_{\lambda, \eta}^2(\x, t) \big( |u_{n} - u^*|^2 + |(u_n - u^*)_t|^2 + |\nabla (u_n - u^*)|^2\big)  d\x dt
	\\
	+  \int_{\Omega_T}	W_{\lambda, \eta}^2(\x, t) \big|c(\x)\mathcal E_{tt} - \Delta \mathcal E \big|^2d\x dt
	 +\epsilon\|u^*\|^2_{H^p_{\lambda, 0}(\Omega_T)}.
	 \label{3.22}
\end{multline}
Recall from \eqref{test fn h} that $z = u_{n + 1} - u^* - \mathcal E$ and the inequality $(a + b)^2 \geq \frac{1}{2}a^2 - b^2$.
Fix $\lambda \geq \frac{\ln C - \ln \varepsilon}{\eta T^2}$ such that $C \lambda^3 e^{-\lambda \eta T^2} < \varepsilon.$
It follows from \eqref{3.22} that
\begin{equation}
	\lambda \|u_{n + 1} - u^*\|^2_{H^1_{\lambda, \eta}(\Omega_T)} 
	\leq C\|u_{n} - u^*\|^2_{H^1_{\lambda, \eta}(\Omega_T)} 
	+ C\lambda \|\mathcal E\|_{H^p_{\lambda, \eta}(\Omega_T)}^2
	 +\epsilon\|u^*\|^2_{H^p_{\lambda, 0}(\Omega_T)}.
\end{equation}
Therefore,
\begin{equation}
	 \|u_{n + 1} - u^*\|^2_{H^1_{\lambda, \eta}(\Omega_T)} 
	\leq \frac{C}{\lambda}\|u_{n} - u^*\|^2_{H^1_{\lambda, \eta}(\Omega_T)} 
	+ C \|\mathcal E\|_{H^p_{\lambda, \eta}(\Omega_T)}^2
	 +\frac{\epsilon}\lambda\|u^*\|^2_{H^p_{\lambda, 0}(\Omega_T)}.
	 \label{3.25}
\end{equation}
Since $\|\mathcal E\|_{H^p(\Omega)} < \delta$,
$\|\mathcal E\|_{H^p_{\lambda, \eta}(\Omega)}^2 \leq e^{2\lambda P_2^2} \delta^2$
where $P_2 = \max_{\x \in \overline \Omega}\{|\x - \x_0|\}$.
 For any $\gamma \in (0, 1)$. If $\delta < e^{-\lambda P_2^2/\gamma}$ then it is obvious that $\|\mathcal E\|_{H^p_{\lambda, \eta}(\Omega_T)}^2 < \delta^{2 - 2\gamma}$. In this case, it follows from \eqref{3.25} that
\begin{equation}
	\|u_{n + 1} - u^*\|^2_{H^1_{\lambda, \eta}(\Omega_T)}  \leq \theta \|u_{n} - u^*\|^2_{H^1_{\lambda, \eta}(\Omega_T)} + C(\delta^{2 - 2\gamma} + \epsilon \|u^*\|_{H^p_{\lambda, 0}}(\Omega_T)^2)
	\label{3.27}
\end{equation}
where $\theta = C/\lambda \in (0, 1).$
Applying \eqref{3.27} with $n$ replaced by $n - 1, \dots, 1$, we obtain \eqref{conv}. The proof is complete.

\section{Numerical study} \label{num}

For simplicity, we numerically solve Problem \ref{ISP} in the 2-dimensional case.
We set $\Omega = (-1, 1)^2$.
To solve the forward problem, we approximate $\R^2$ by the square $G = (-R, R)^2$ where $R > 1$ is a large number. 
We then employ the explicit scheme to solve \eqref{main eqn} on $G$. 
In our computation, $R = 4$.
We set the final time $T = 2 < 4$. 
In this setting, the nonlinear wave generated by an unknown source supported inside $\Omega$ does not have enough time to propagate to $\partial G.$ 
Therefore, the computed wave  is the restriction of the solution to \eqref{main eqn}  on $G \times [0, T]$.
Our implementation is based on the finite difference method, in which the step size in space is $d_\x = 0.03$ and the step size in time is $d_t = 0.0132.$ These values of step sizes in space and time are chosen such that $d_t/d_\x = 0.44 < 0.5$. Thus, the computed wave by the explicit scheme is reliable. 

In this section, we set the noise level $\delta = 10\%.$ That means the data we use are given by 
\begin{equation}
	f(\x, t) = f^*(\x, t)(1 + \delta {\rm rand})
	\quad
	\mbox{and}
	\quad
	g(\x, t) = g^*(\x, t)(1 + \delta {\rm rand})
	\label{noisenoise}
\end{equation}
for all $(\x, t) \in \Gamma_T$ and rand is the function taking uniformly distributed random number in the range $[-1, 1].$
In all of our numerical examples, the stopping number $n_0$ in Step \ref{n0} of Algorithm \ref{alg} is set to be 7. It is not necessary to set $n_0$ as a bigger number. This is because of the fast convergence of the exponential rate $O(\theta^n)$ as $n \to \infty$ for some $\theta \in (0, 1)$. See Theorem \ref{thm} and \eqref{conv}.

Implementing Algorithm \ref{alg} is straight forward. 
The most complicated part is Step \ref{step update}. That is how to minimize a function $J_n$ defined in \eqref{Jn}.
This step is equivalent to the problem solving a linear equation subject to some boundary and initial constraints by the quasi-reversibility method. 
We employ the optimization package built-in in MATLAB with the command ``lsqlin" for this step. 
The implementation is very similar to the one described in \cite[Section 6]{LocNguyen:ip2019}.
We do not repeat the details of implementation here.
\begin{remark}
Although the norm in the regularization term is $H^p_{\lambda, 0}(\Omega_T)$, in our computational program, we use a simpler norm $H^2_{\lambda, 0}(\Omega_T)$ to simplify the computational scripts. 
  This change does not effect the satisfactory outputs. 
Moreover, this change significantly simplifies the implementation.
If one uses the finite element method, (s)he does not have to use high degree for the finite element space.
  The "artificial" parameters for all tests below are: $\epsilon = 10^{-14}$, $\lambda = 2.1$ and $\eta = 0.5$. These parameters are numerically chosen such that the reconstruction result is satisfactory for test 1. Then, we use these values for all other tests.
 \end{remark}

Regarding the initial solution $u_0$, we have mentioned in Step \ref{s1} of Algorithm \ref{alg} that $u_0$ can be any function in $H$. That means we can take any function in $H^p(\Omega_T)$ such that $u|_{\Gamma_T} = f$, $\partial_{\nu} u|_{\Gamma_T}= g$ and $u_t|_{\Omega \times \{t = 0\}} = 0.$
A natural way to compute such a function $u_0$ is to find a function $u$ satisfying
\begin{equation}
	\left\{
		\begin{array}{rcll}
			\Delta u(\x, t)&=& 0 &(\x, t) \in \Omega_T,\\
			u(\x, t) &=& f &(\x, t) \in \Gamma_T,\\
			u(\x, t) &=& g &(\x, t) \in \Gamma_T,\\
			u_t(\x, 0) &=& 0 &\x \in \Omega.
		\end{array}
	\right.
\label{u0}
\end{equation}
We do so by the quasi-reversibility method. That means we set $u_0$ as the minimizer of
\[
	\int_{\Omega_T} W^2_{\lambda, \eta} (\x, t)|\Delta u|^2 dt d\x + \epsilon \|u\|_{H^2_\lambda(\Omega_T)}
\]  
subject to the constraints as in the last three lines of  \eqref{u0}.
Again, we refer the reader to \cite[Section 6]{LocNguyen:ip2019} for the details in implementation of the quasi-reversibility method using the optimization package of MATLAB.

We show four (4) numerical tests. 

\begin{figure}[h!]
	\subfloat[The true source function $p^*$.]{\includegraphics[width=.3\textwidth]{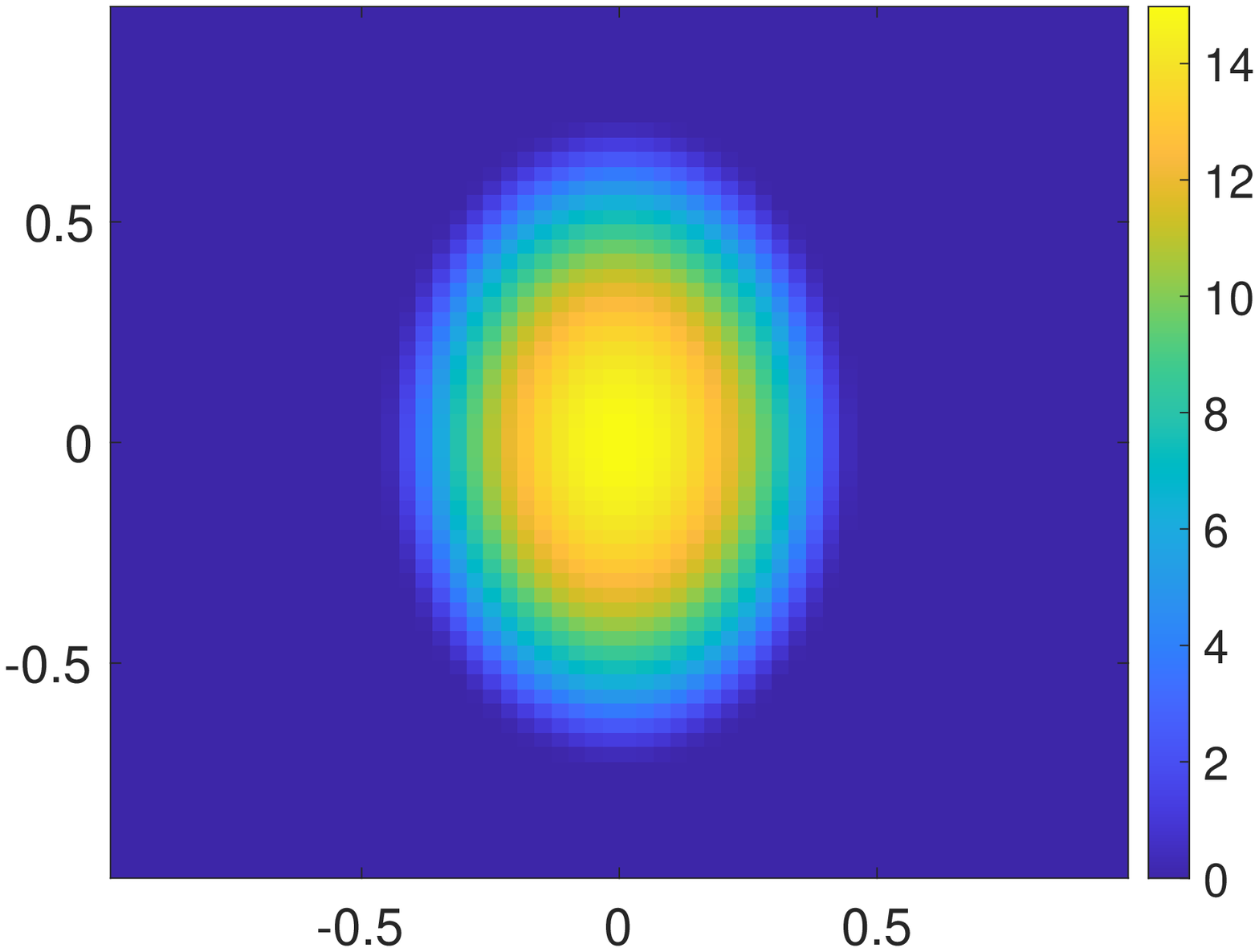}}
	\quad
	\subfloat[\label{fig 1b} The initial solution $p_{\rm init} = u_0(\cdot, 0)$ computed by solving \eqref{u0} by the quasi-reversibility method.]{\includegraphics[width=.3\textwidth]{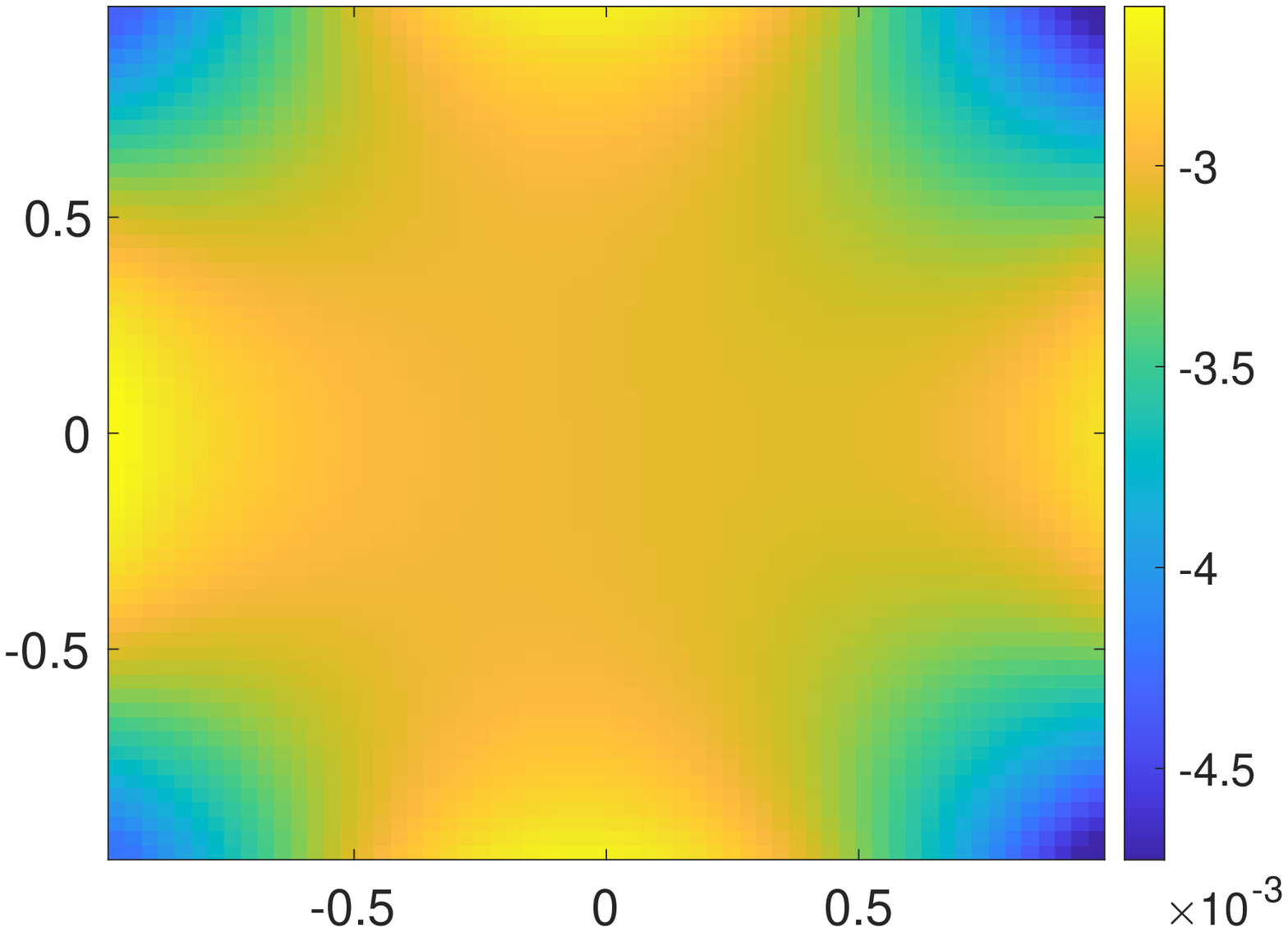}}
	
	\subfloat[\label{fig 1c}The computed source function $p_{\rm comp}$.]{\includegraphics[width=.3\textwidth]{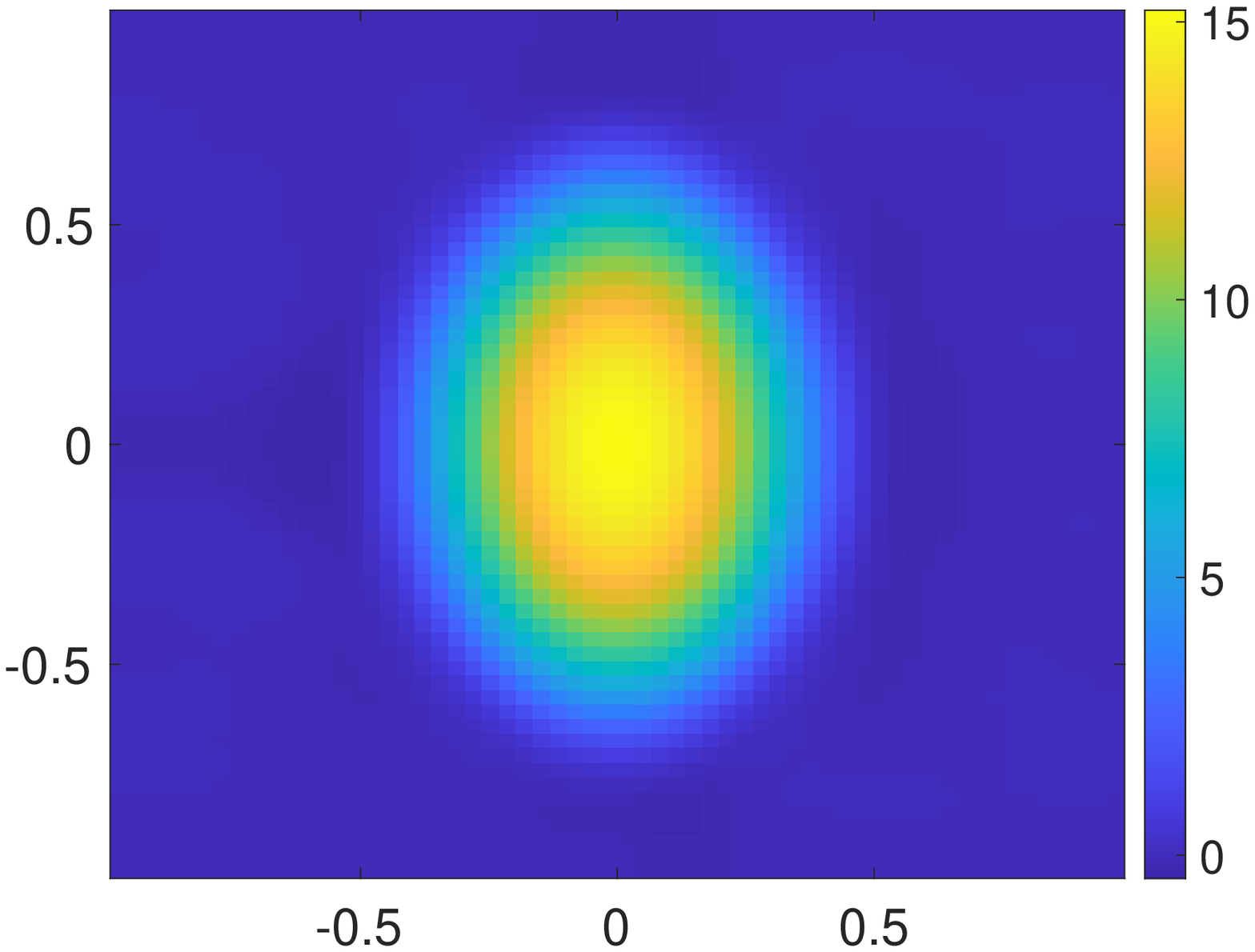}}
	\quad
	\subfloat[\label{fig 1d}The relative difference $\frac{|p_{\rm comp} - p^*|}{\|p^*\|_{L^{\infty}(\Omega)}}$.]{\includegraphics[width=.3\textwidth]{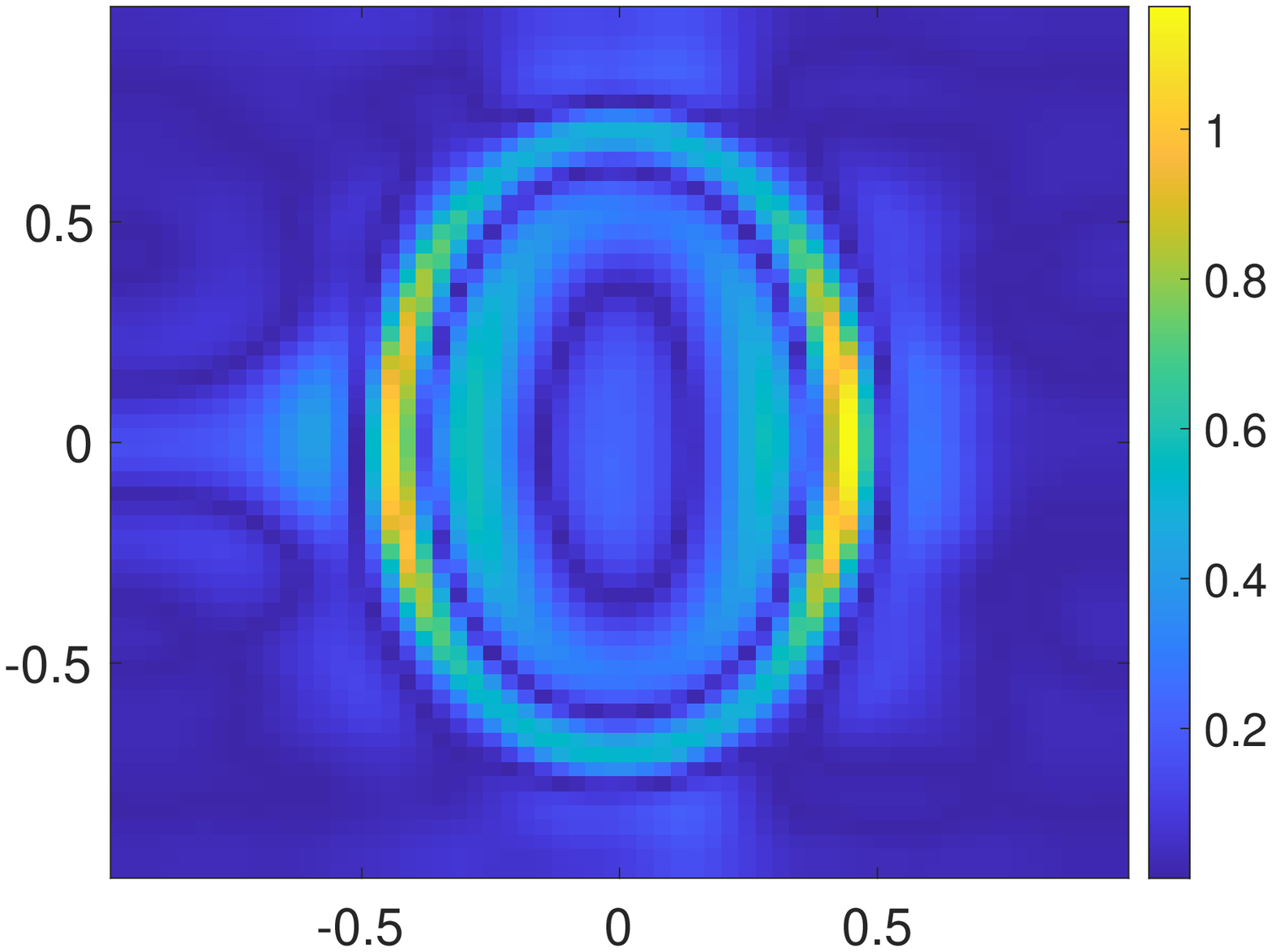}}
	\quad
	\subfloat[\label{fig 1e}The consecutive difference $\|p_{n} - p_{n-1}\|_{L^{\infty}(\Omega)}$.]{\includegraphics[width=.3\textwidth]{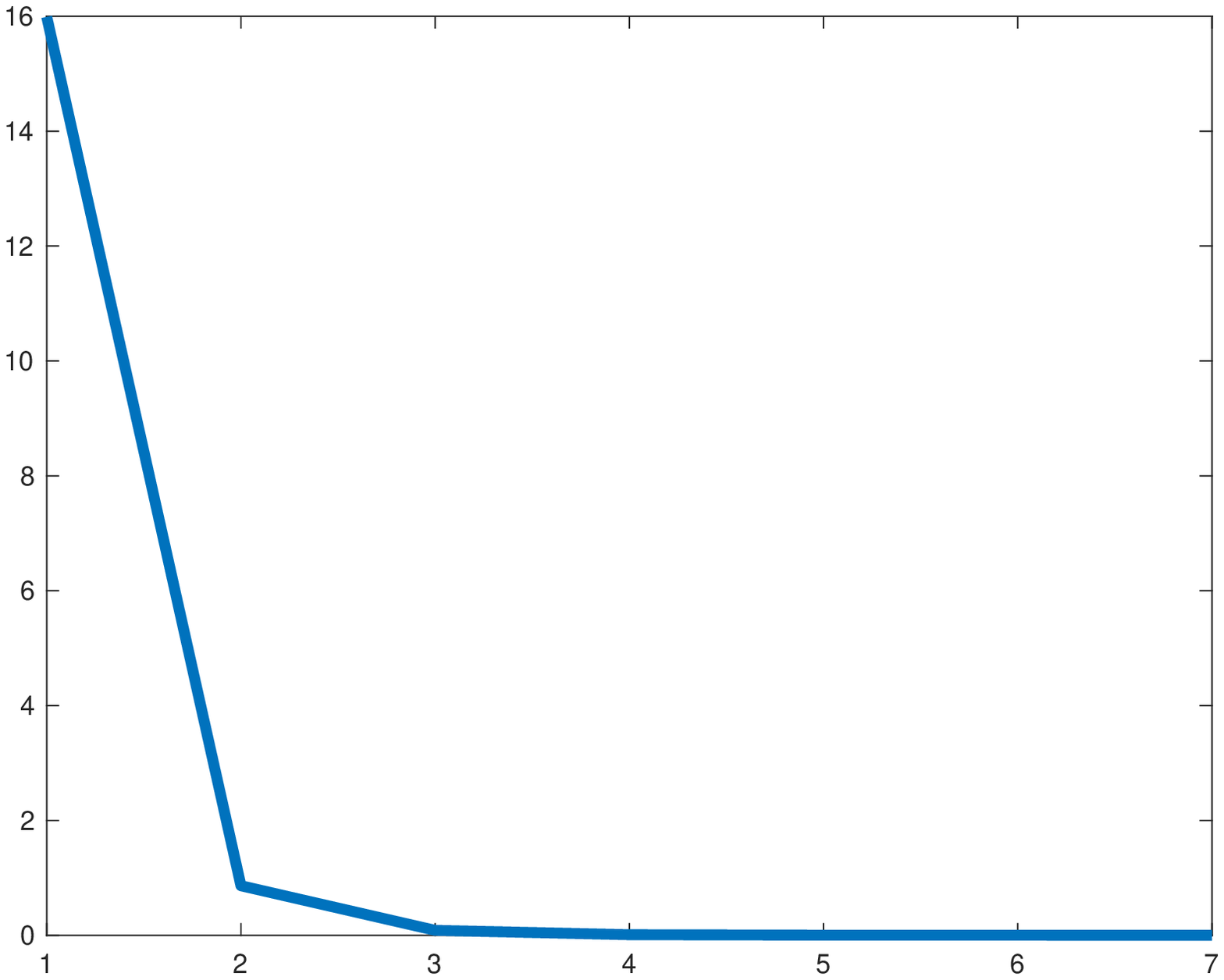}}
	\caption{\label{fig 1} Test 1. The true and computed source functions. 
	The computed source function can be successfully found although the support of the source function is large and the value of this function is 15, which is high.
	The initial solution $p_{\rm init}$ is faraway from $p^*$.}
\end{figure}

{\bf Test 1.} 
In this test, we set the nonlinearity $f$ to be
\[	
	f(\x, t, u, u_t, \nabla u) = \sqrt{u^2 + 1} + |\nabla u|
\]
and the true source function is
\[
	p^*(x, y) =  \left\{
		\begin{array}{ll}
			15 e^{\frac{r^2(x, y)}{r^2(x, y) - 1}} &\mbox{if} \,  r(x, y)< 1,\\
			0 &\mbox{otherwise}
		\end{array}
	\right.
\]
where $r(x, y) = \frac{x^2}{0.5^2} + \frac{y^2}{0.8^2}$. 
The support of the true source function $p^*$ is an ellipse centered at the origin with the $x-$radius $0.5$ and $y-$radius 0.8. The maximum value of $p^*$ is 15. 
In this test, the true function $p^*$ is far away from the background. So, finding an initial guess is impossible.
The numerical result of this test is displayed in Figure \ref{fig 1}.

It is evident from Figures \ref{fig 1c}--\ref{fig 1d} that we successfully compute the source function without requesting an initial guess.
In fact, the initial  solution $p_{init}$ obtained by solving \eqref{u0} by the quasi-reversibility method is almost equal to the zero function. It is far away from the true solution $p^*$, see Figure \ref{fig 1b}.
It follows from Figure \ref{fig 1e} that our method converges fast after only 4 iterations. This observation is consistent to the conclusion of Theorem \ref{thm}.
 The $L^2$-relative error $\frac{\|p_{\rm comp} - p^*\|_{L^2(\Omega)}}{\|p^*\|_{L^2(\Omega)}}$ in Test 1 is 5.4\%.

\begin{figure}[h!]
	\subfloat[The true source function $p^*$.]{\includegraphics[width=.3\textwidth]{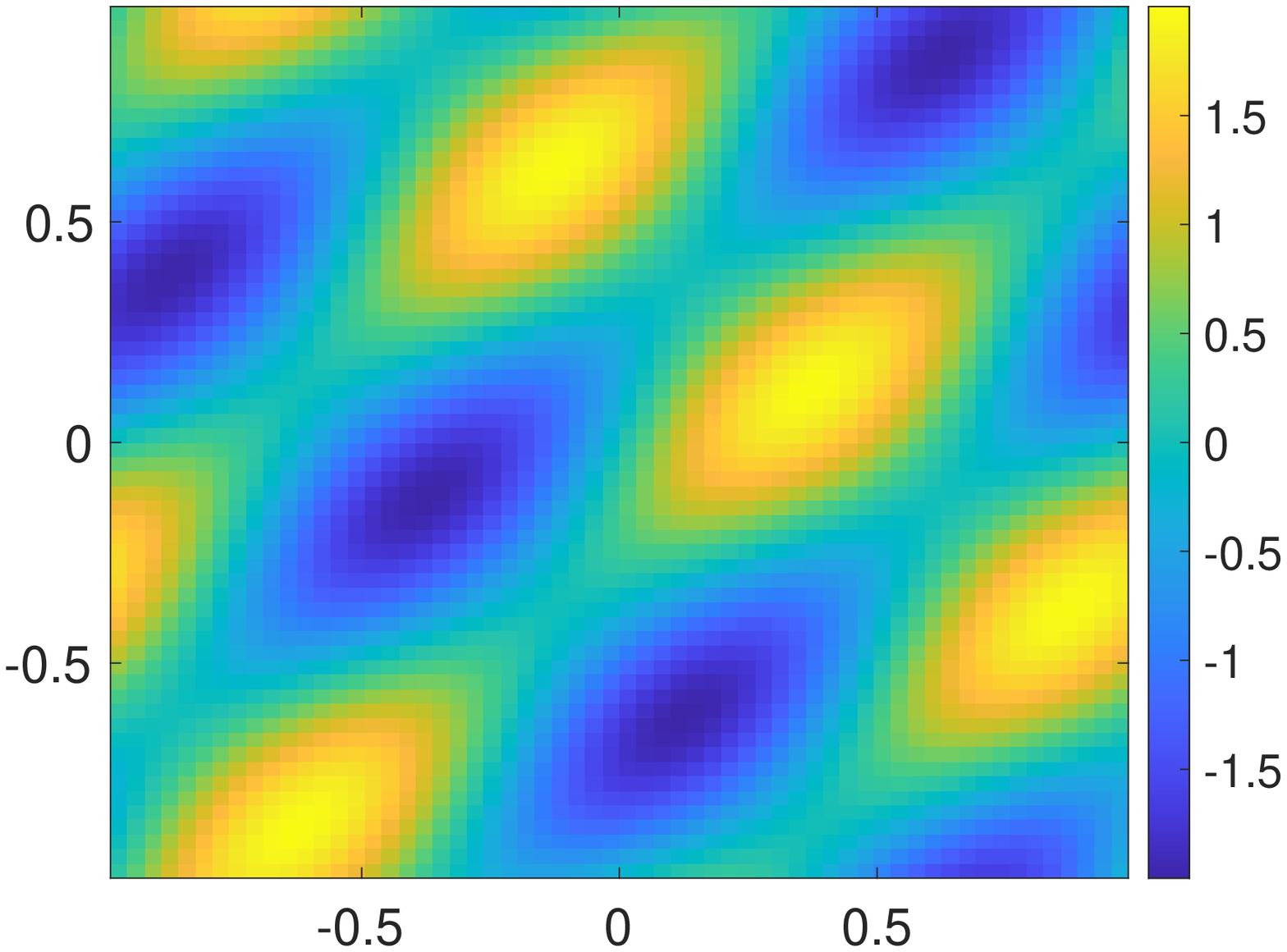}}
	\quad
	\subfloat[\label{fig 2b} The initial solution $p_{\rm init} = u_0(\cdot, 0)$ computed by solving \eqref{u0} by the quasi-reversibility method.]{\includegraphics[width=.3\textwidth]{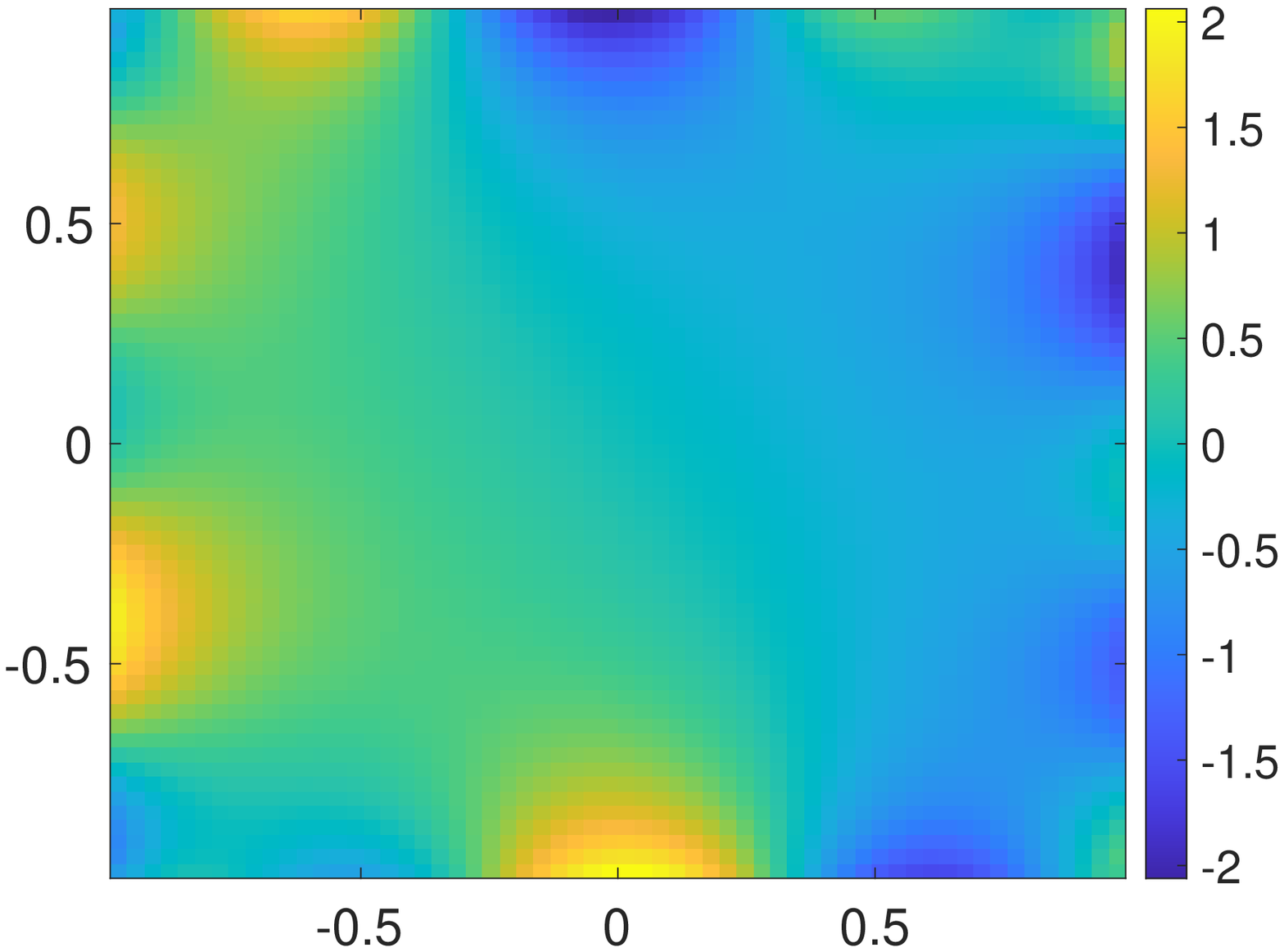}}
	
	\subfloat[\label{fig 2c}The computed source function $p_{\rm comp}$.]{\includegraphics[width=.3\textwidth]{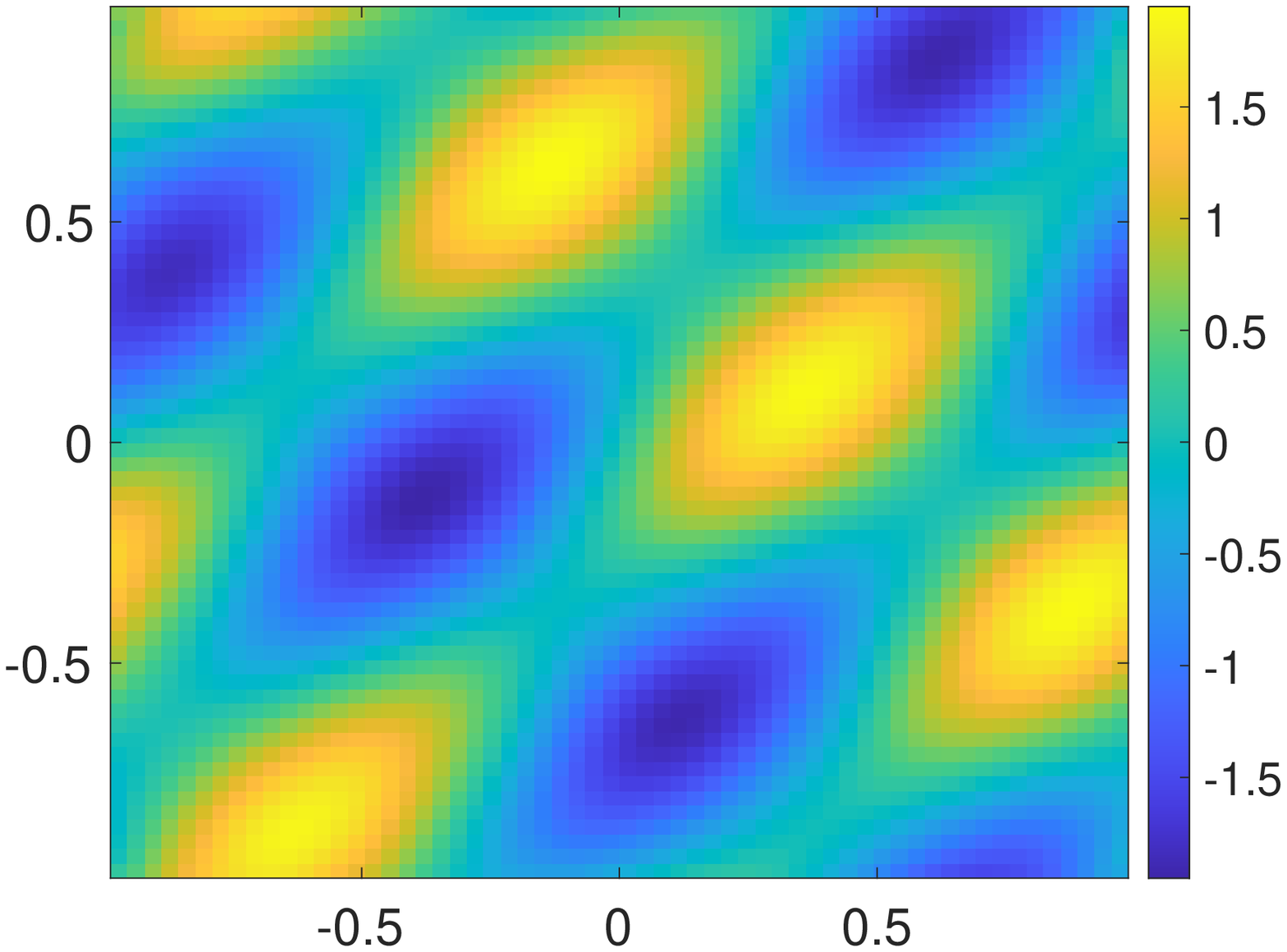}}
	\quad
	\subfloat[\label{fig 2d}The relative difference $\frac{|p_{\rm comp} - p^*|}{\|p^*\|_{L^{\infty}(\Omega)}}$.]{\includegraphics[width=.3\textwidth]{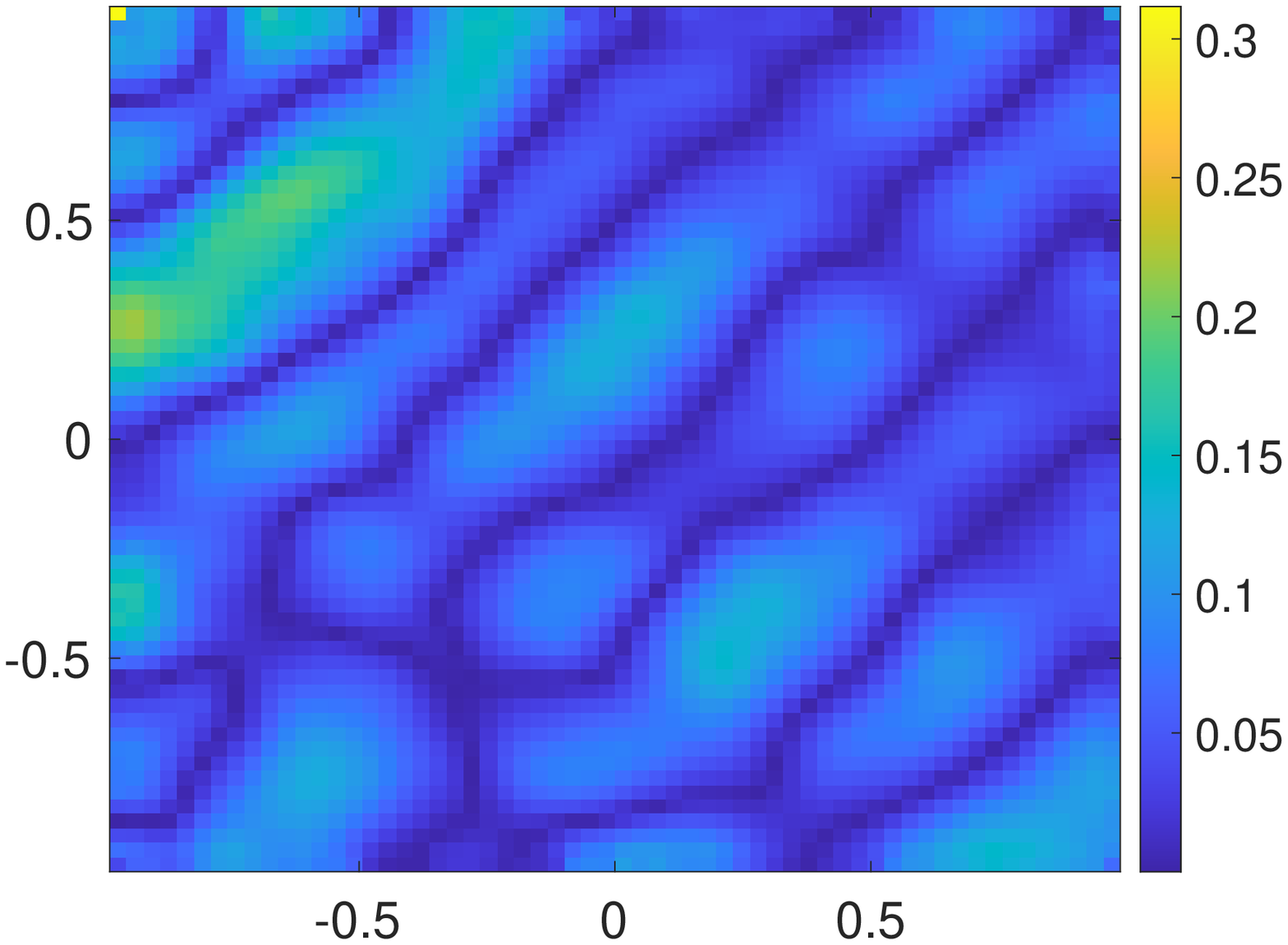}}
	\quad
	\subfloat[\label{fig 2e}The consecutive difference $\|p_{n} - p_{n-1}\|_{L^{\infty}(\Omega)}$.]{\includegraphics[width=.3\textwidth]{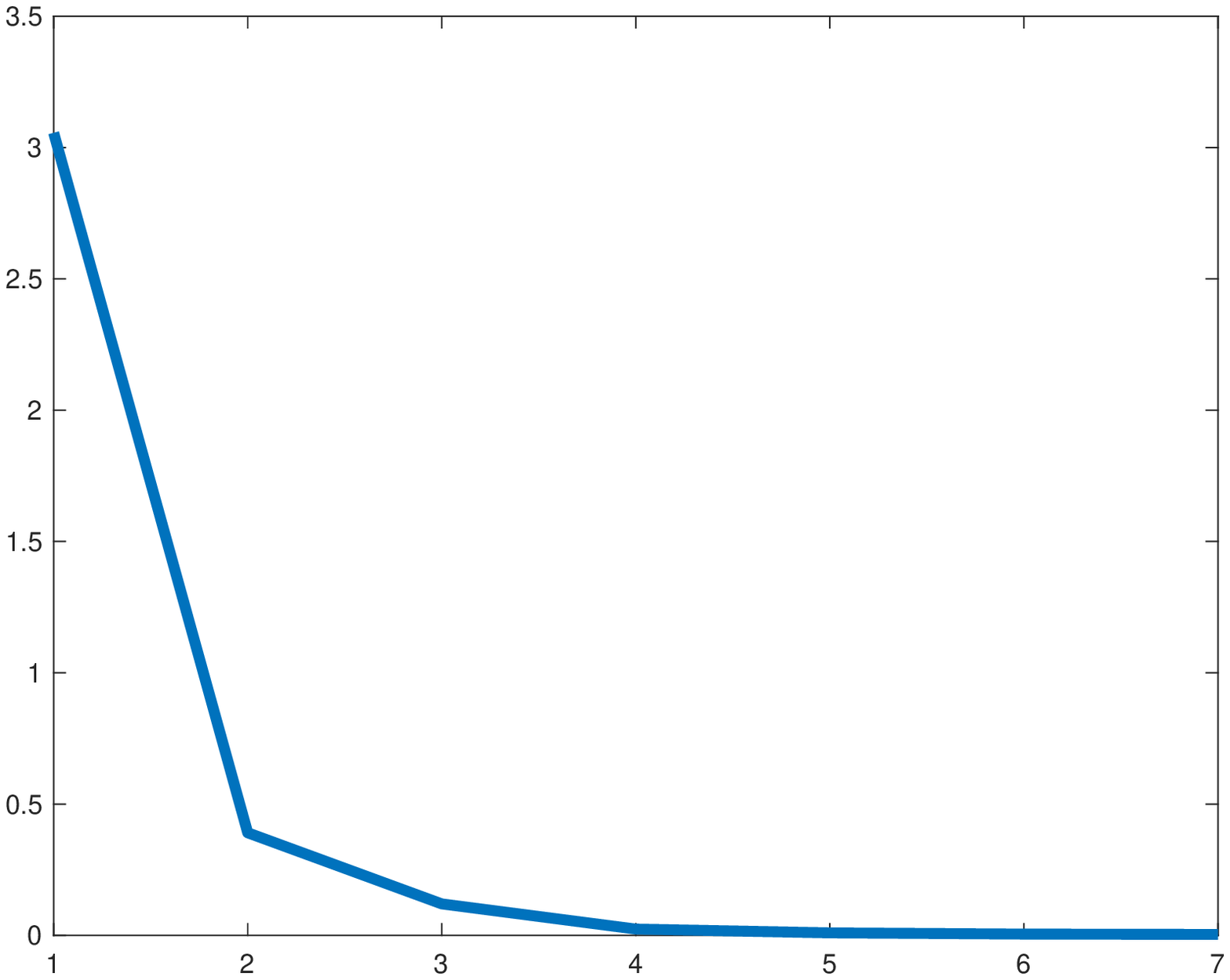}}
	\caption{\label{fig 2} Test 2. The true and computed source functions. 
	The computed source function can be successfully found although the geometry structure of $p^*$ is complicated and the nonlinearity $f$ is not in the class $C^1$.}
\end{figure}

 {\bf Test 2.} We consider the source function with more complicated structure. 
 The true source function is given by 
 \[
 	p^*(x, y) = \sin(\pi(x + y)) +  \sin((2\pi(x-y))) \quad
	\mbox{for all } \x = (x, y) \in \Omega.
 \]
 The nonlinearity in this test is a bounded function
 \[
 	f(\x, t, u, u_t, \nabla u) = \min\{e^u + |\nabla u|, 10\}.
 \]
Although the nonlinearity $f$ is not in the class $C^1$ and the geometric structure of the true source function is complicated, one will see that Algorithm \ref{alg} is effective, see Figure \ref{fig 2}.

As in Test 1, it is evident from Figures \ref{fig 2c}--\ref{fig 2d} that we successfully compute the source function without requesting an initial guess.
In fact, the initial  solution $p_{init}$ obtained by solving \eqref{u0} by the quasi-reversibility method is  far away from the true solution $p^*$, see Figure \ref{fig 2b}.
It follows from Figure \ref{fig 2e} that our method converges fast after only 5 iterations. This observation is consistent to the conclusion of Theorem \ref{thm}.
 The $L^2$-relative error $\frac{\|p_{\rm comp} - p^*\|_{L^2(\Omega)}}{\|p^*\|_{L^2(\Omega)}}$ in Test 1 is 7.07\%.

\begin{figure}[h!]
	\subfloat[The true source function $p^*$.]{\includegraphics[width=.3\textwidth]{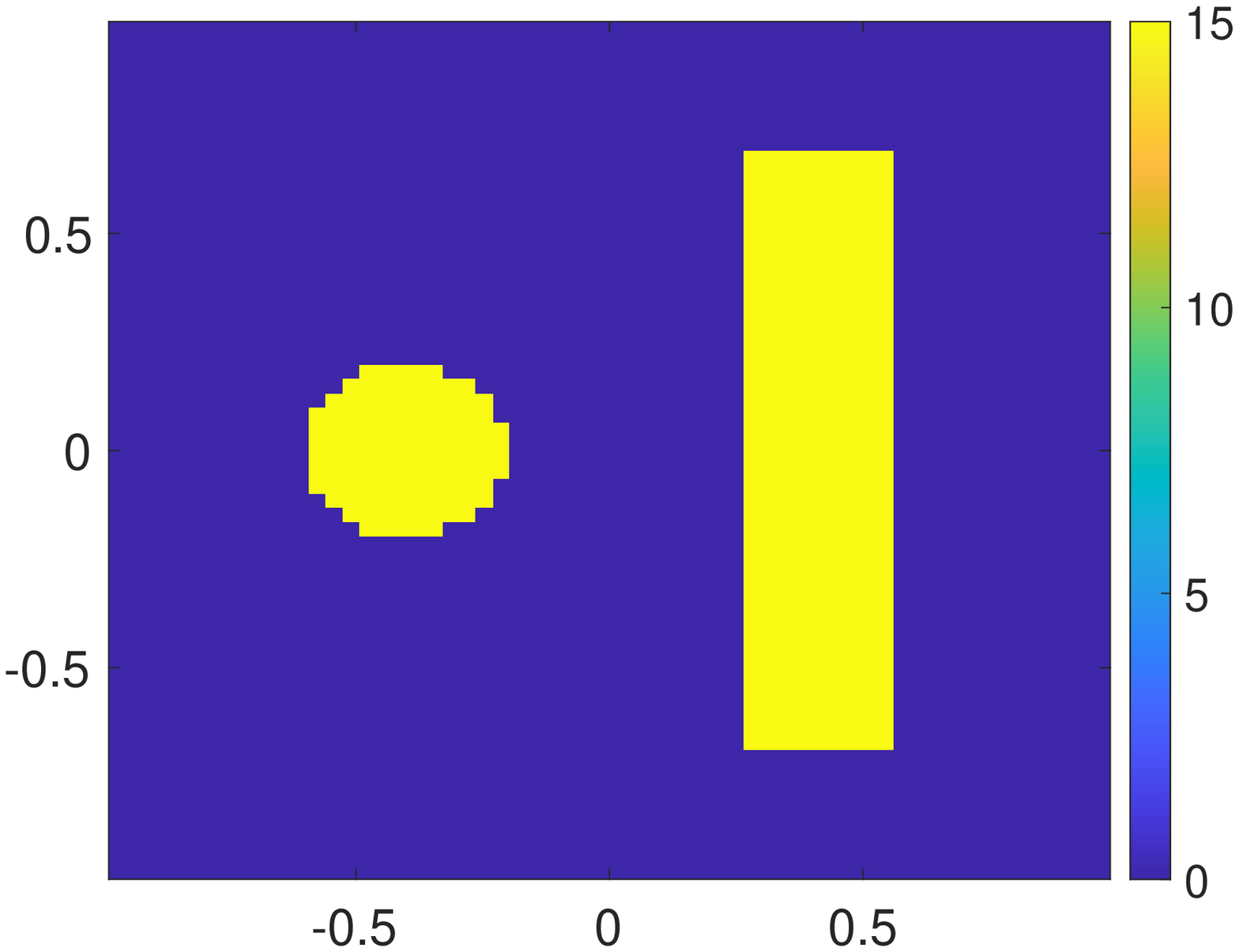}}
	\quad
	\subfloat[\label{fig 3b} The initial solution $p_{\rm init} = u_0(\cdot, 0)$ computed by solving \eqref{u0} by the quasi-reversibility method.]{\includegraphics[width=.3\textwidth]{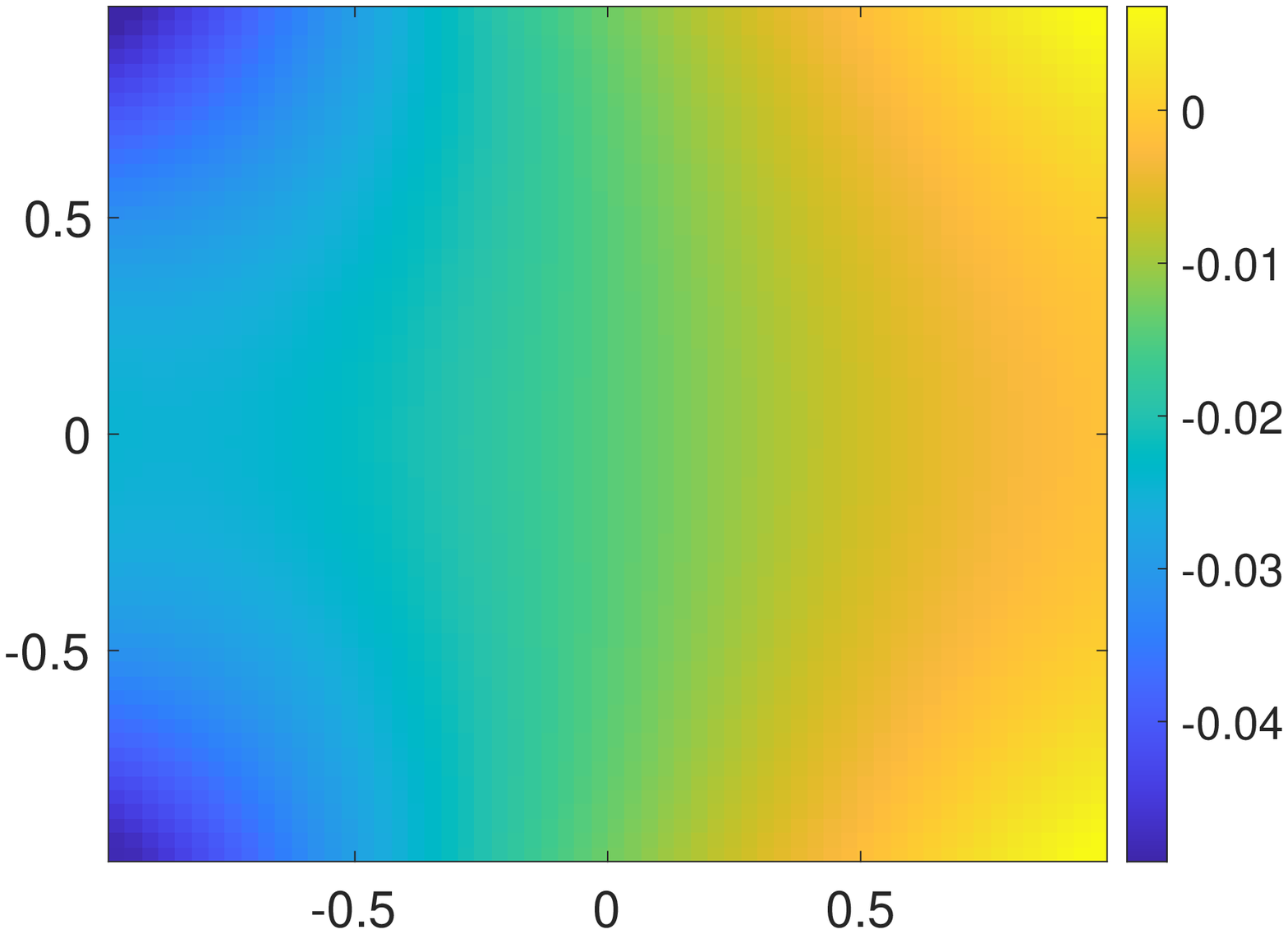}}
	
	\subfloat[\label{fig 3c}The computed source function $p_{\rm comp}$.]{\includegraphics[width=.3\textwidth]{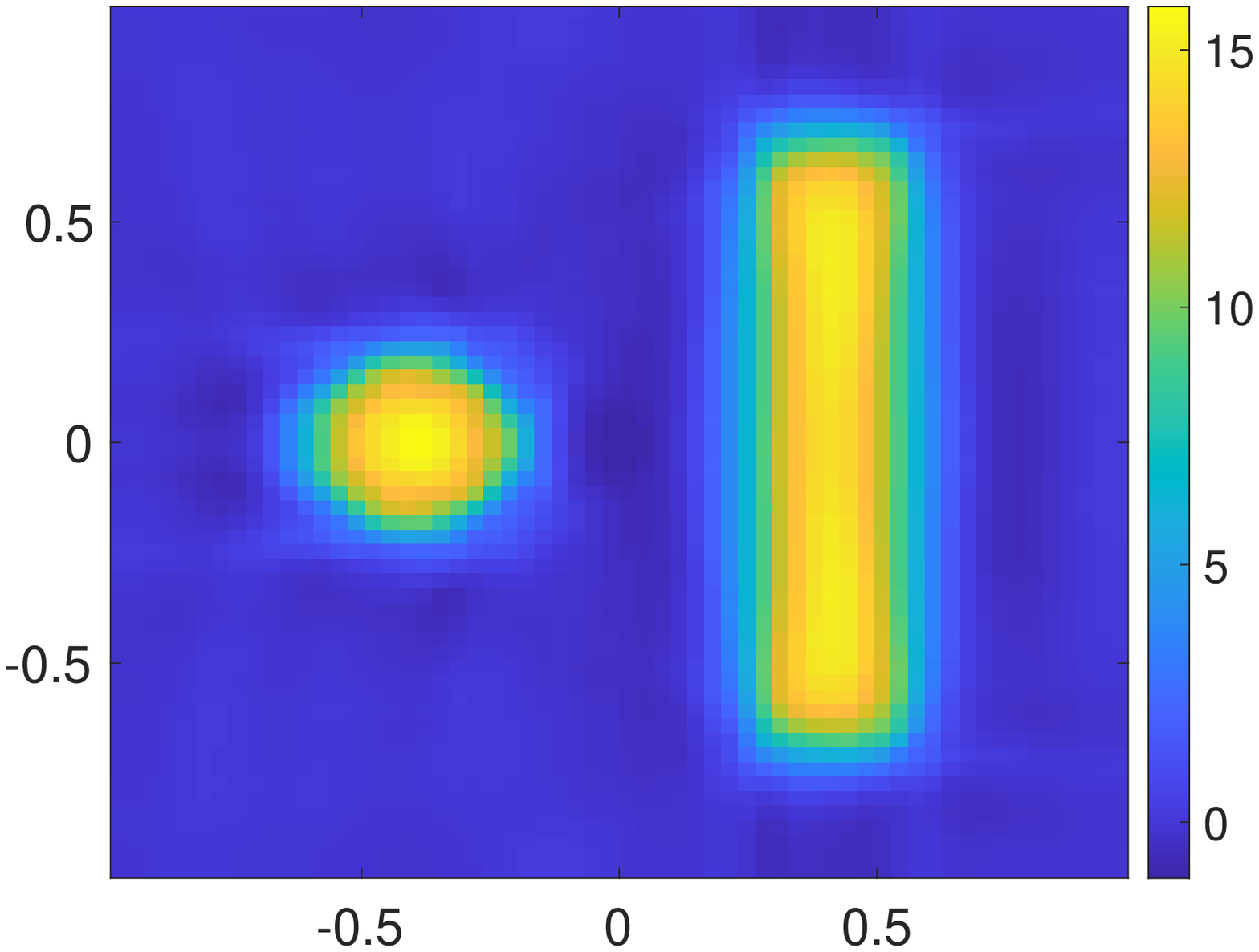}}
	\quad
	\subfloat[\label{fig 3d}The relative difference $\frac{|p_{\rm comp} - p^*|}{\|p^*\|_{L^{\infty}(\Omega)}}$.]{\includegraphics[width=.3\textwidth]{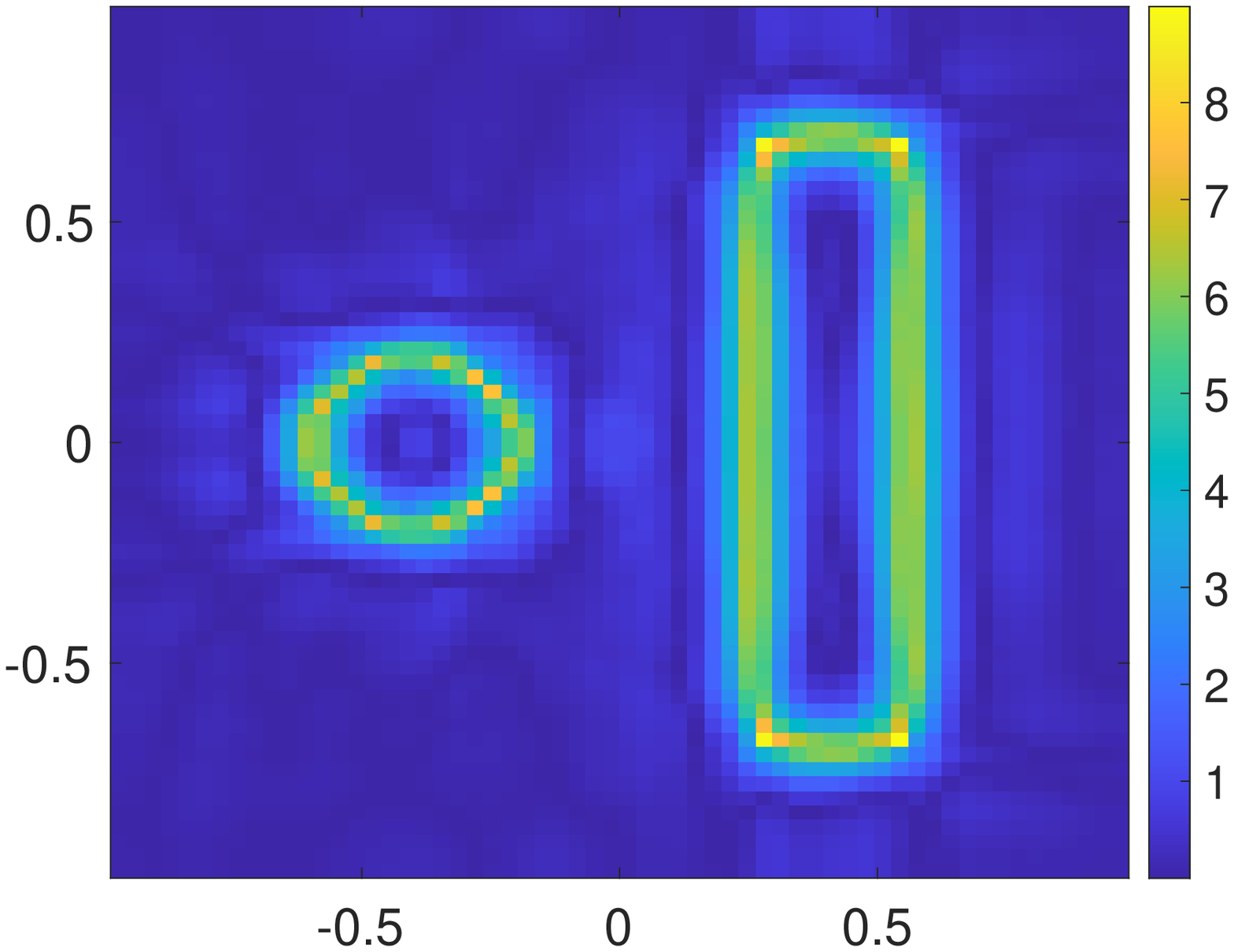}}
	\quad
	\subfloat[\label{fig 3e}The consecutive difference $\|p_{n} - p_{n-1}\|_{L^{\infty}(\Omega)}$.]{\includegraphics[width=.3\textwidth]{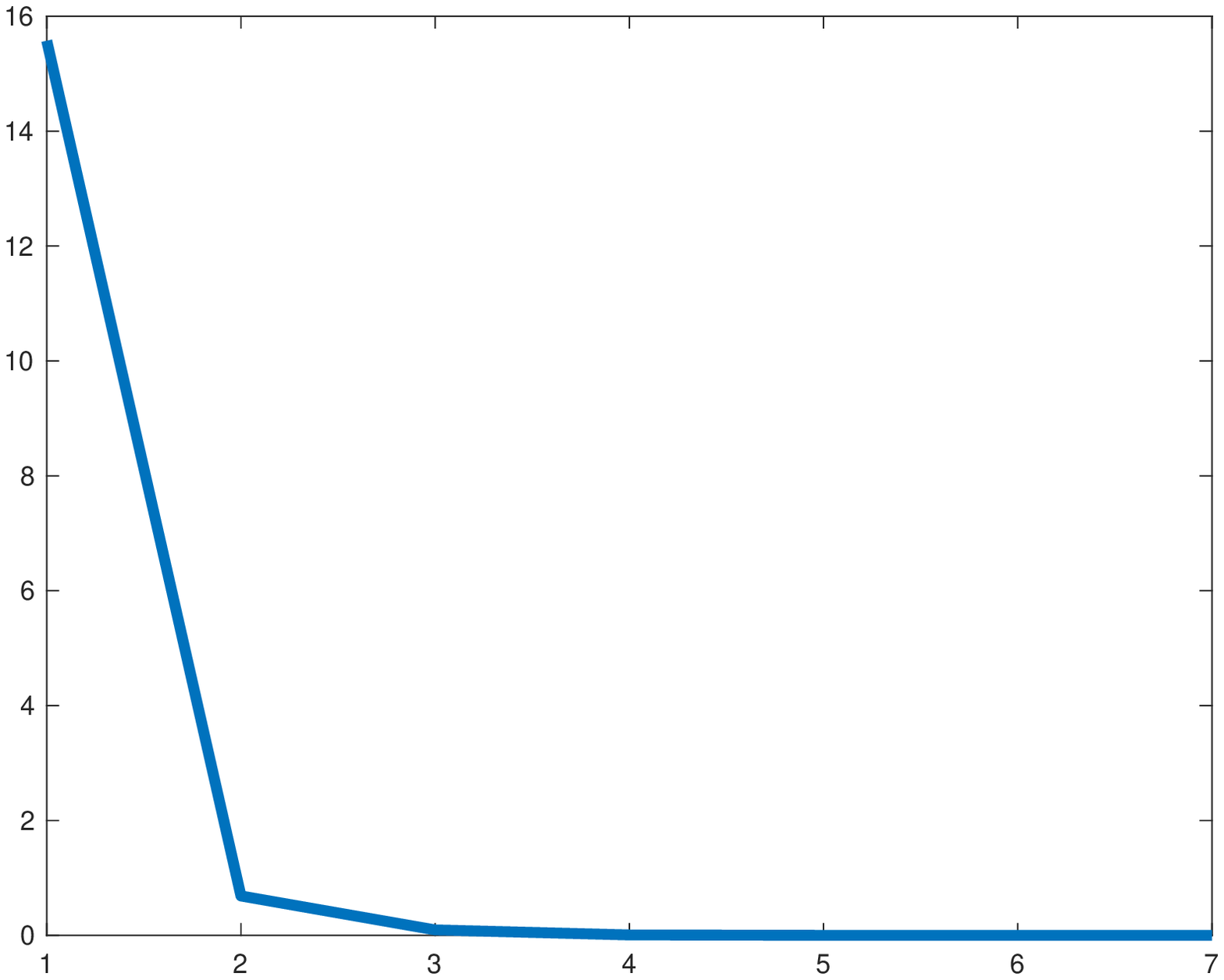}}
	\caption{\label{fig 3} Test 3. The true and computed source functions. 
	The computed source function can be successfully found although $p^*$ is far away from the background zeros function and the nonlinearity $f$ is not in the class $C^1$ and not bounded.}
\end{figure}

{\bf Test 3.}
We test the case when the true source function is not smooth, given by the following piece-wise constant function
\[
	p^*(x, y) = 
	\left\{
	\begin{array}{ll}
		15 &\mbox{if} \max\Big\{\frac{|x - 0.4|}{0.15}, \frac{|y|}{0.7} \Big\} < 1 \, \mbox{or } (x + 0.4)^2 + y^2 < 0.04,\\
		0 &\mbox{otherwise}. 
	\end{array}
	\right.
\]
The support of $p^*$ involves a vertical rectangle and a disk. The maximum value of $p^*$ is 15. Finding an initial guess for this function is impossible.
In this test, we use the nonlinearity 
\[
	f(\x, t, u, u_t, \nabla u) = \min\{u^2 + 1,10\} + |\nabla u|.
\]
Like Test 2, the nonlinearity in this test is not in the class $C^1$ but unlike Test 2, $f$ is unbounded.
The computed solution to Problem \ref{ISP} of Test 3 is displayed in Figure \ref{fig 3}.

The numerical result is satisfactory. One can see that we can effectively reconstruct both rectangular and circular ``inclusions". 
The maximum value of $p_{\rm comp}$ is 15.8 (relative error is $5.3\%$).
Unlike the case when $p^*$ is smooth, the $L^2$ relative error in this test is  high 34.5\%. However, 
one can see from Figure \ref{fig 3d} that the large error occurs significantly on the borders of the inclusions while the error is small almost everywhere. 
Like in Test 1 and Test 2, the convergence rate of Algorithm \ref{alg} is fast.
One can see from Figure \ref{fig 3e} that we can obtain the computed source function with only 4 iterations.

\begin{figure}[h!]
	\subfloat[The true source function $p^*$.]{\includegraphics[width=.3\textwidth]{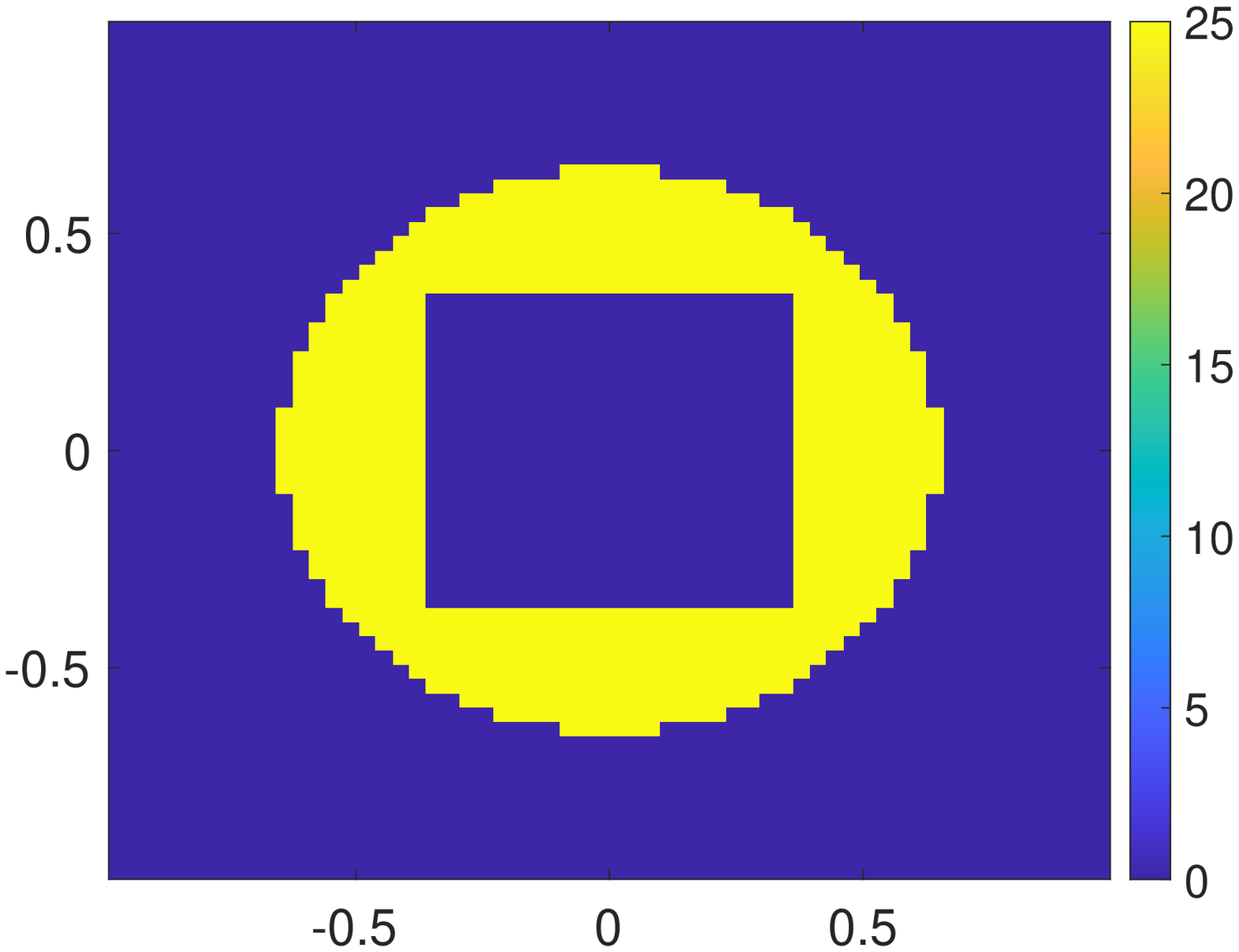}}
	\quad
	\subfloat[\label{fig 3b} The initial solution $p_{\rm init} = u_0(\cdot, 0)$ computed by solving \eqref{u0} by the quasi-reversibility method.]{\includegraphics[width=.3\textwidth]{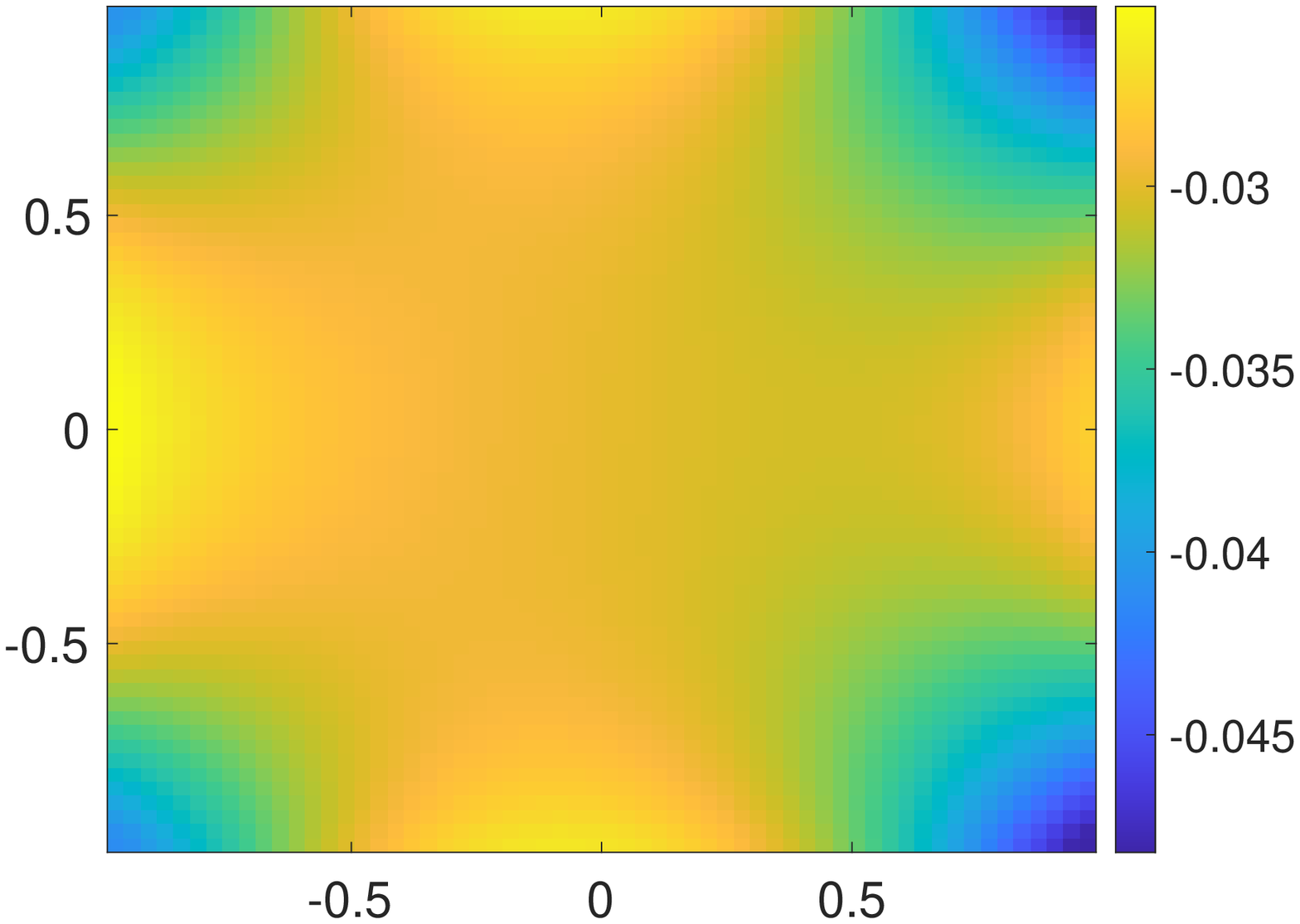}}
	
	\subfloat[\label{fig 4c}The computed source function $p_{\rm comp}$.]{\includegraphics[width=.3\textwidth]{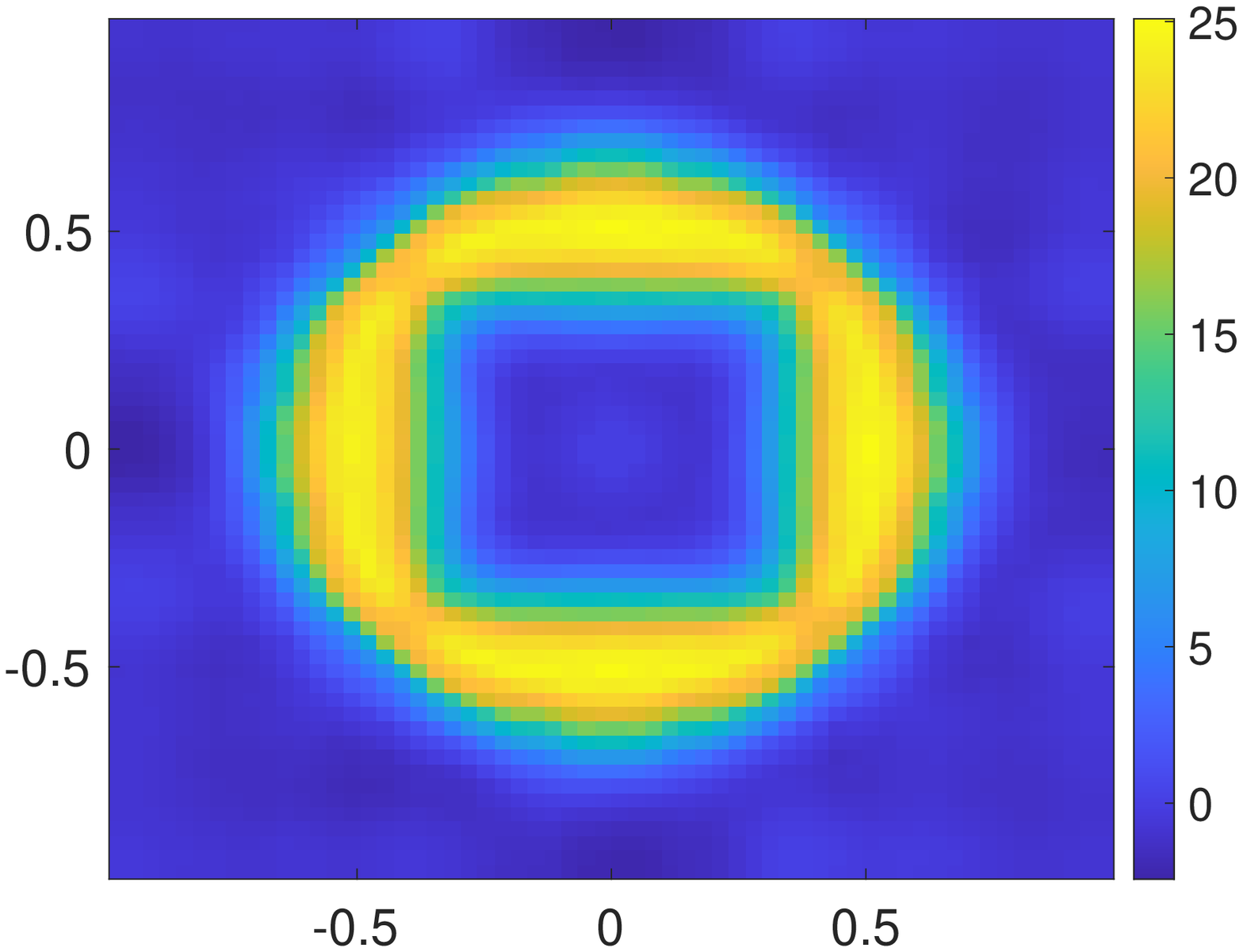}}
	\quad
	\subfloat[\label{fig 4d}The relative difference $\frac{|p_{\rm comp} - p^*|}{\|p^*\|_{L^{\infty}(\Omega)}}$.]{\includegraphics[width=.3\textwidth]{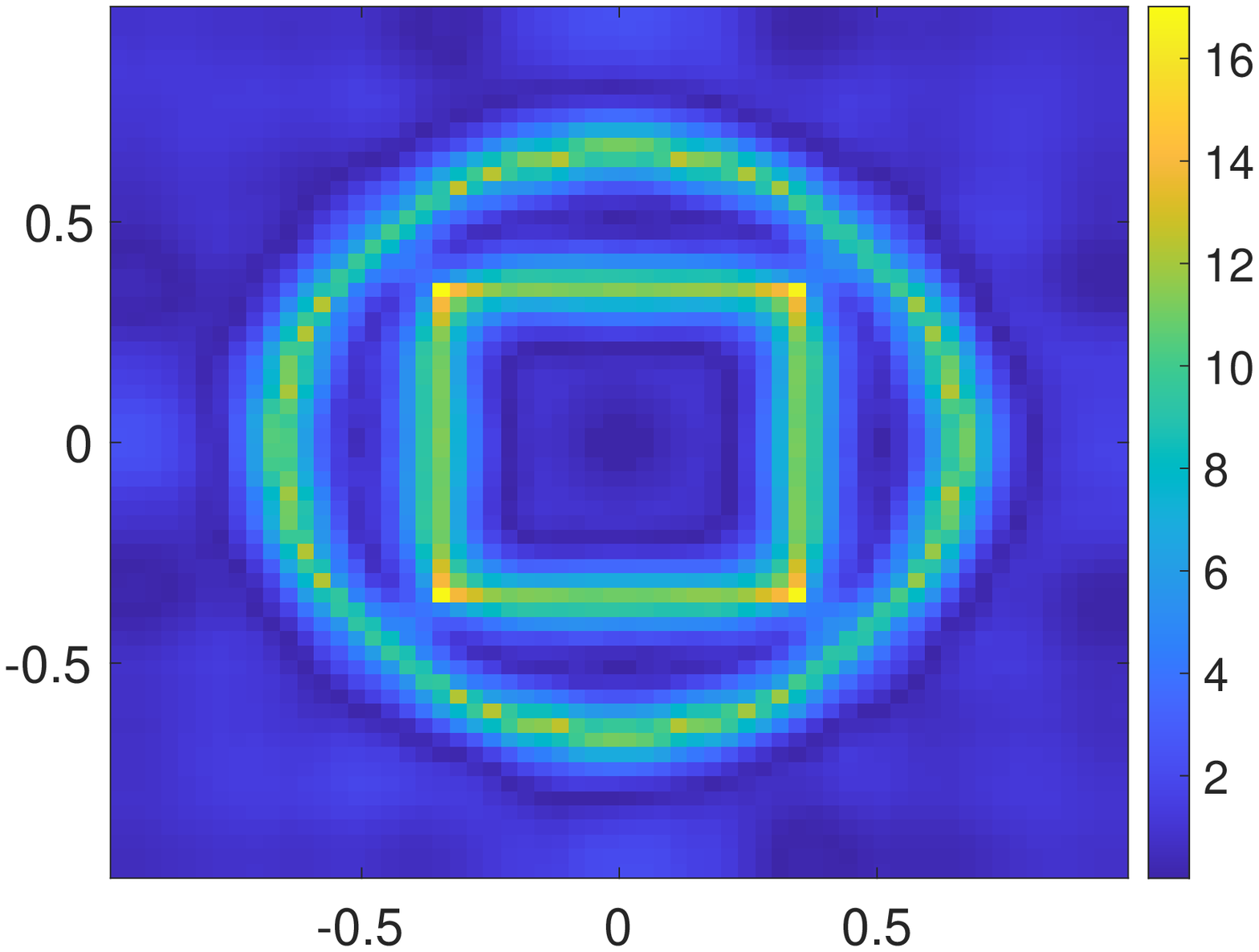}}
	\quad
	\subfloat[\label{fig 4e}The consecutive difference $\|p_{n} - p_{n-1}\|_{L^{\infty}(\Omega)}$.]{\includegraphics[width=.3\textwidth]{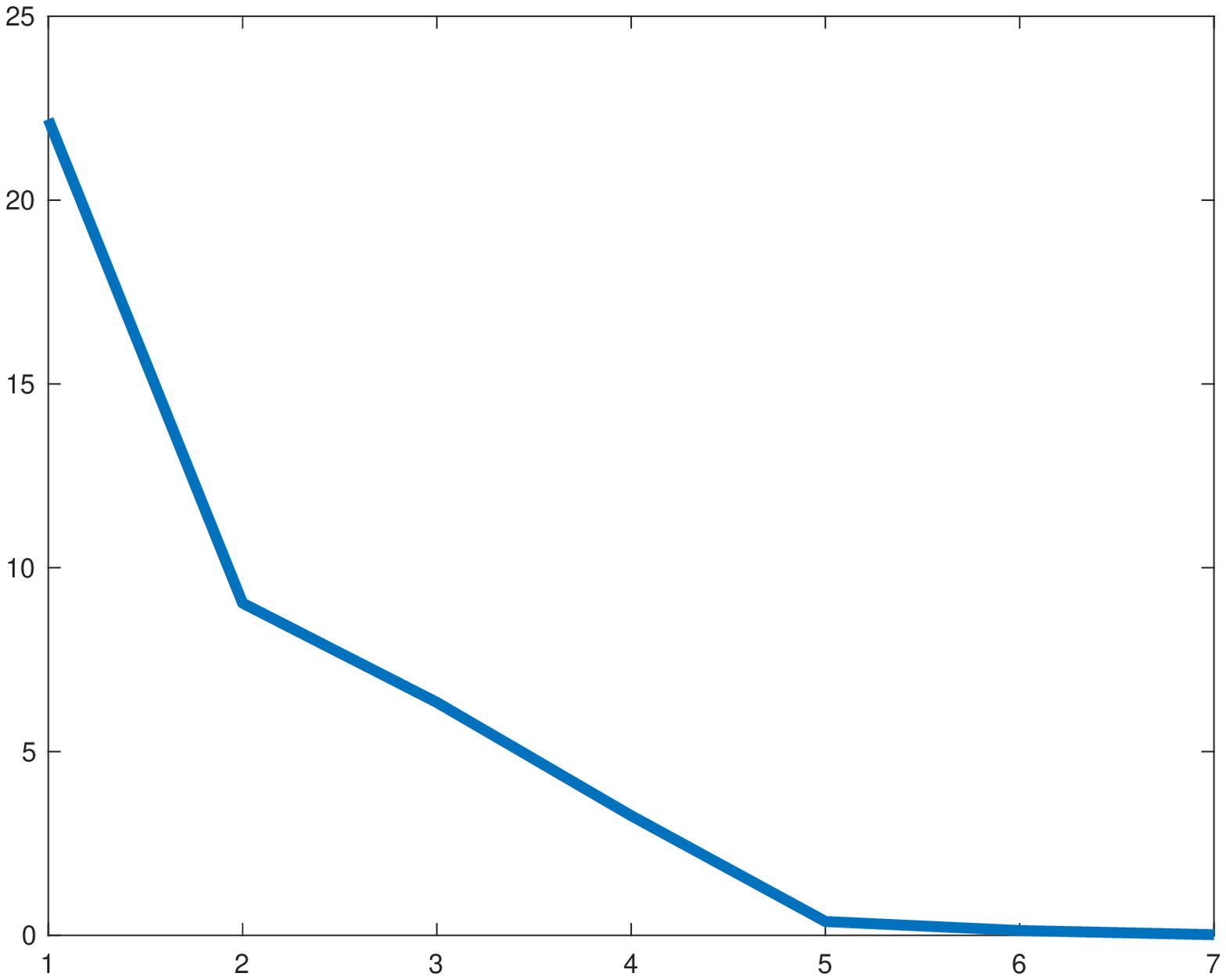}}
	\caption{\label{fig 4} Test 4. The true and computed source functions. 
	The computed source function can be successfully found although $p^*$ is far away from the background zeros function and the nonlinearity $f$ is not in the class $C^1$ and not bounded.}
\end{figure}

{\bf Test 4.}
We now consider a more complicated case when the true source function has large contrast and its  circular support has a rectangular void. 
The function $p^*$ is given by
\[
	p^*(x, y) = 
	\left\{
		\begin{array}{ll}
			25 & \mbox{if}\, x^2 + y^2 < 0.65^2 \mbox{and} \, \max\{|x|, |y|\} > 0.35,\\
			0 & \mbox{otherwise}.
		\end{array}
	\right.
\]
The nonlinearity in this test is set to be
\[
	f(\x, t, u, u_t, \nabla u) = \min\{e^u, 100\} + |\nabla u|.
\]
The nonlinearity in this test is not bounded nor in the class $C^1$. 
The numerical result for this test is in Figure \ref{fig 4}.

Although the true source function $p^*$ and the nonlinearity in this test is complicated, the obtained numerical result is out of expectation. One can see that we can effectively reconstruct the shape of the support of the source function. 
The maximum value of $p_{\rm comp}$ is 25.087 (relative error is $0.4\%$).
Unlike the case when $p^*$ is smooth, the $L^2$ relative error in this test is  high 36.0\%. However,
one can see from Figure \ref{fig 4d} that the large error occurs significantly on the borders of the support of $p^*$ while the error is small almost everywhere. 
Like in all tests above, the convergence rate of Algorithm \ref{alg} is fast.
One can see from Figure \ref{fig 4e} that we can obtain the computed source function with only 7 iterations.

\section{Concluding remarks} \label{sec 5}

We have solve an inverse problem that is to reconstruct the initial condition of nonlinear hyperbolic equation. 
Although this problem is nonlinear, we solve it globally. That means we do not require a good initial guess.
We first define an operator $\Phi$ such that the true solution to the inverse problem is the fixed point of $\Phi$.
We construct a recursive sequence $\{u_n\}_{n \geq 0}$ whose initial term $u_0$ can be taken arbitrary and the $n^{\rm th}$ term $u_n = \Phi(u_{n - 1})$.
We next apply a Carleman estimate to prove the convergence of this sequence.
The convergence is of the exponential rate.
Moreover, we have proved that the stability of our method with respect to noise is of the H\"older type.
Numerical results are satisfactory. 
\section*{Acknowledgement} This work  was supported by US Army Research Laboratory and US Army Research
Office grant W911NF-19-1-0044.


\end{document}